\documentclass [11pt] {article}
\usepackage{amssymb}
\title{Universality of Fedosov's Construction for\\
       Star Products of Wick Type on\\
       Pseudo-K\"ahler Manilfolds}

\author{{\bf Nikolai Neumaier
\thanks{Nikolai.Neumaier@physik.uni-freiburg.de}}\\[3mm]
Fakult\"at f\"ur Mathematik und Physik\\Universit\"at
Freiburg\\Hermann-Herder-Stra\ss e 3\\D-79104 Freiburg i.~Br.\\
Germany\\[3mm]}
\date{April 2002\\[3mm]FR-THEP-2002/07}
\newcommand {\BEQ} [1] {\begin {equation} \label{#1}}
\newcommand {\EEQ} {\end {equation}}
\newcommand {\BEQAR} [1] {\begin {eqnarray} \label{#1}}
\newcommand {\EEQAR} {\end {eqnarray}}
\newcommand {\TinyWi} {{\mbox{\rm \tiny Wick}}}
\newcommand {\TinyAW} {{\overline{\TinyWi}}}
\newcommand {\TinyF} {{\mbox{\rm \tiny F}}}
\newcommand {\TinyK} {{\mbox{\rm \tiny K}}}

\newcommand {\tauWi} {\tau_\TinyWi}
\newcommand {\tauAW} {\tau_\TinyAW}
\newcommand {\tauK} {\tau_\TinyK}
\newcommand {\starWi} {\mathbin{\star_{\TinyWi}}}
\newcommand {\starAW} {\mathbin{\star_{\TinyAW}}}
\newcommand {\circF} {\mathbin{\circ_{\TinyF}}}
\newcommand {\circWi} {\mathbin{\circ_{\TinyWi}}}
\newcommand {\circAW} {\mathbin{\circ_{\TinyAW}}}
\newcommand {\circK} {\mathbin{\circ_{\TinyK}}}

\newcommand {\rWi} {r_{\TinyWi}}
\newcommand {\rAW} {r_{\TinyAW}}
\newcommand {\rK} {r_{\TinyK}}
\newcommand {\DF} {\mathfrak D_{\TinyF}}
\newcommand {\DWi} {{\mathfrak D_{\TinyWi}}}
\newcommand {\DAW} {\mathfrak D_{\TinyAW}}
\newcommand {\DK} {\mathfrak D_{\TinyK}}
\newcommand {\astF} {\mathbin{\ast_{\TinyF}}}
\newcommand {\astWi} {\mathbin{\ast_{\TinyWi}}}
\newcommand {\astAW} {\mathbin{\ast_{\TinyAW}}}
\newcommand {\astK} {\mathbin{\ast_{\TinyK}}}
\newcommand {\OF} {{\Omega_{\TinyF}}}
\newcommand {\OWi} {{\Omega_{\TinyWi}}}
\newcommand {\OAW} {{\Omega_{\TinyAW}}}
\newcommand {\sF} {{s_{\TinyF}}}
\newcommand {\sWi} {{s_{\TinyWi}}}
\newcommand {\sAW} {{s_{\TinyAW}}}
\newcommand {\cc} [1] {\overline {{#1}}}
\newcommand {\Hol} [1] {\mathcal O ({#1})}
\newcommand {\AHol} [1] {\cc{\mathcal O} ({#1})}
\newcommand {\Cinf} [1] {\mathcal C^\infty ({#1})}
\newcommand {\Ginf} [1] {\Gamma^\infty ({#1})}
\newcommand {\im} {{\mathrm{i}}}
\newcommand {\id} {{\mathrm{id}}}
\newcommand {\Lie} {{\mathcal L}}
\newcommand {\ad} {{\mathrm{ad}}}
\newcommand {\WL} {\mbox{$\mathcal W \! \otimes \! \Lambda$}}
\newcommand {\deltaH} {\delta_{\mbox{\rm \tiny H}}}
\newcommand {\adF} {{\ad_{\TinyF}}}
\newcommand {\adWi} {{\ad_{\TinyWi}}}
\newcommand {\adAW} {{\ad_{\TinyAW}}}
\newcommand {\adK} {{\ad_{\TinyK}}}
\newcommand {\Deg} {\mathrm{Deg}}
\newcommand {\C} {{\mathrm C}}
\newcommand {\Pa} {{\mathrm P}}
\newcommand {\pa} {{\mbox{\rm\tiny P}}}
\newcommand {\opp} {{\mbox{\rm\tiny opp}}}
\newcommand {\starfl} {\mathbin{\star_{\mbox{\rm\tiny opp}}}}
\renewcommand {\d} {\mathrm{d}}
\newenvironment {PROOF}{\small {\sc Proof:}}{{\hspace*{\fill}
                       $\square$}}
\newenvironment {INNERPROOF}{\small {\sc Proof:}}{{\hspace*{\fill}
                            $\bigtriangledown$}}
\newtheorem {LEMMA} {Lemma} [section]
\newtheorem {SUBLEMMA} [LEMMA] {Sublemma}
\newtheorem {PROPOSITION} [LEMMA] {Proposition}
\newtheorem {THEOREM} [LEMMA] {Theorem}
\newtheorem {COROLLARY} [LEMMA] {Corollary}
\newtheorem {DEFINITION} [LEMMA] {Definition}
\newtheorem {REMARK} [LEMMA] {Remark}
\newtheorem {DEDUCTION} [LEMMA] {Deduction}
\begin{document}
\maketitle
\begin{abstract}
In this paper we construct star products on a pseudo-K\"ahler
manifold $(M,\omega,I)$ using a modification of the Fedosov method
based on a different fibrewise product similar to the Wick product
on $\mathbb C^n$. Having fixed the used connection to be the
pseudo-K\"ahler connection these star products shall depend on
certain data given by a formal series of closed two-forms on $M$
and a certain formal series of symmetric contravariant tensor
fields on $M$. In a first step we show that this construction is
rich enough to obtain star products of every equivalence class by
computing Deligne's characteristic class of these products. Among
these products we uniquely characterize the ones which have the
additional property to be of Wick type which means that the
bidifferential operators describing the star products only
differentiate with respect to holomorphic directions in the first
argument and with respect to anti-holomorphic directions in the
second argument. These star products are in fact strongly related
to star products with separation of variables introduced and
studied by Karabegov. This characterization gives rise to special
conditions on the data that enter the Fedosov procedure. Moreover,
we compare our results that are based on an obviously coordinate
independent construction to those of Karabegov that were obtained
by local considerations and give an independent proof of the fact
that star products of Wick type are in bijection to formal series
of closed two-forms of type $(1,1)$ on $M$. Using this result we
finally succeed in showing that the given Fedosov construction is
universal in the sense that it yields all star products of Wick
type on a pseudo-K\"ahler manifold. Due to this result we can make
some interesting observations concerning these star products; we
can show that all these star products are of Vey type and in
addition we can uniquely characterize the ones that have the
complex conjugation incorporated as an anti-automorphism.
\end{abstract}
\newpage
\tableofcontents
\section{Introduction}
\label{IntroSec}
The concept of deformation quantization as introduced in the
pioneering articles \cite{BayFla78} by Bayen, Flato, Fr\o nsdal,
Lichnerowicz and Sternheimer has proved to be an extremely useful
framework for the problem of quantization: the question of
existence of star products $\star$ (i.e. formal, associative
deformations of the classical Poisson algebra of complex-valued
functions $\Cinf{M}$ on a symplectic or more generally, on a
Poisson manifold $M$, such that in the first order of the formal
parameter $\nu$ the commutator of the star product yields the
Poisson bracket) has been answered positively by DeWilde and
Lecomte \cite{DeWLec83b}, Fedosov \cite{Fed94}, Omori, Maeda and
Yoshioka \cite{OmoMaeYos91} in the case of a symplectic phase space
as well as by Kontsevich \cite{Kon97} in the more general case of a
Poisson manifold. Moreover, star products have been classified up
to equivalence in terms of geometrical data of the phase space by
Nest and Tsygan \cite{NT95a}, Bertelson, Cahen and Gutt
\cite{BerCahGut97}, Weinstein and Xu \cite{WX97} on symplectic
manifolds and the classification on Poisson manifolds is due to
Kontsevich \cite{Kon97}. Comparisons between the different results
on classification and reviews can be found in articles of Deligne
\cite{Del95}, Gutt and Rawnsley \cite{GR99,Gut00}, Neumaier
\cite{Neu99} and the thesis of Halbout \cite{Hal99a}.

In the case of deformation quantizations with separation of
variables on K\"ahler manifolds Karabegov proved existence and gave
a classification using a formal deformation of the K\"ahler form in
\cite{Kar96,Kar98}. Moreover, he has shown in \cite{Kar99} that the
Fedosov approach to star products of Wick type that only differ
from those with separation of variables by interchanging the
r\^{o}les of holomorphic and anti-holomorphic directions and a
minus sign in front of the symplectic resp. K\"ahler form
considered by Bordemann and Waldmann in \cite{BorWal96a}
corresponds to the trivial deformation in the sense of his
classification. This result naturally raises the question, whether
this construction can be generalized as to obtain all star products
of Wick type on a K\"ahler manifold.

The aim of this paper is to give a positive answer to this question
stressing the power and beauty of Fedosov's method and allowing for
a global, geometrical description of all star products of Wick type
that avoids local considerations in coordinates. In addition the
further investigation of additional properties of the star products
of Wick type is extremely simplified due to the simple, lucid
description of these star products encoded in the Fedosov
derivation. Moreover, our result represents the justification to
restrict to the Fedosov framework when addressing further questions
like the construction of representations for star products of Wick
type.

The paper is organized as follows: In Section \ref{NotDefSec} we
establish our notation and remember some basic definitions that are
needed in the context of star products of Wick type on
pseudo-K\"ahler manifolds. Section \ref{FedProdEquCharClaSec} is
devoted to the presentation of the Fedosov constructions using
different fibrewise products. Moreover, we give explicit
constructions for equivalence transformations relating the
resulting star products that enable us using our result in
\cite[Thm. 4.4]{Neu99} to show that the given construction with
$\circWi$ is rich enough to obtain star products in every
equivalence class. After this first important result we give a
unique characterization of the star products of Wick type in
Section \ref{WickCharSec} obtained from our construction in terms
of the data that enter the Fedosov construction. In addition we can
prove that the star products of Wick type only depend on a formal
series of closed two-forms of type $(1,1)$. In Section
\ref{KaraCharSec} we remember some results of Karabegov obtained in
\cite{Kar96} and give elementary proofs of the important
statements. Using these results we are in the position to prove the
main result of the paper in Section \ref{FedUnivSec} stating that
every star product of Wick type (not only such a star product in
every possible equivalence class) can be obtained by some adapted
Fedosov construction. Using this fact we can draw some additional
conclusions concerning the structure of such products, in
particular we show that each such product is of Vey type. In
Appendix \ref{ComStrConSec} we have collected some notations and
results concerning the splittings into holomorphic and
anti-holomorphic part of the mappings that are essential for the
Fedosov construction that make use of the presence of the complex
structure $I$ and that are important for the proofs and statements
given in Sections \ref{WickCharSec} and \ref{FedUnivSec}. A further
Appendix \ref{BanFixPktSec} is added for completeness providing
some details on the application of Banach's fixed point theorem in
the framework of the Fedosov construction which constitutes an
important and simplifying tool for several proofs in Sections
\ref{FedProdEquCharClaSec} and \ref{WickCharSec}.

{\bf Conventions:} By $\Cinf{M}$ we denote the complex-valued
smooth functions and similarly $\Ginf{T^*M}$ stands for the
complex-valued smooth one-forms et cetera. Moreover, we use
Einstein's summation convention in local expressions.

{\bf Acknowledgements:} I would like to express my thanks to Stefan
Waldmann for a careful reading of the manuscript and many comments
and useful discussions. Moreover, I should like to thank Alexander
V. Karabegov for the inspiring discussions that initiated the
investigations of the topic of the paper. Finally, I want to thank
the DFG-Graduiertenkolleg ``Nichtlineare
Differentialglei\-chun\-gen -- Modellierung, Theorie, Numerik,
Visualisierung'' for financial support.
\section{Some Notations and Basic Definitions}
\label{NotDefSec}
Let $(M,\omega,I)$ denote a pseudo-K\"ahler manifold with
$\dim_{\mathbb R}(M)=2n$, i.e. $(M,g)$ is a pseudo-Riemannian (no
positivity is required) manifold such that the almost complex
structure $I$ is an isometry with respect to $g$ ($g$ is Hermitian)
and the almost complex structure is flat with respect to the
Levi-Civita connection $\nabla$ corresponding to $g$. Under these
conditions it is known that $M$ is a complex manifold such that the
almost complex structure $I$ coincides with the canonical complex
structure of $M$ and that $\omega$ defined by $\omega(X,Y):= g
(IX,Y)$ for $X,Y \in\Ginf{TM}$ is a closed non-degenerate two-form
and hence a symplectic form on $M$, the so-called pseudo-K\"ahler
form. Moreover, one obviously has that $\nabla$ defines a torsion
free symplectic connection since $\omega$ is also flat with respect
to the connection $\nabla$ which is called the pseudo-K\"ahler
connection. As $g$ is Hermitian both $g$ and $\omega$ are of type
$(1,1)$ and in a local holomorphic chart of $M$ one can write
$\omega=\frac{\im}{2}g_{k\cc{l}}\d z^k\wedge \d\cc{z}^l$ with a
Hermitian non-degenerate matrix $g_{k\cc{l}}=2g(Z_k,\cc{Z}_l)$.
Here $Z_k=\partial_{z^k}$ and $\cc{Z}_l=\partial_{\cc{z}^l}$ denote
local base vector fields of type $(1,0)$ and of type $(0,1)$ (cf.
Appendix \ref{ComStrConSec}) that locally span the $+\im$ and
$-\im$ eigenspaces $TM^{1,0}$ and $TM^{0,1}$ of the complex
structure $I$.

The most simple example of a (pseudo-)K\"ahler manifold is given by
$\mathbb C^n$ endowed with the canonical (pseudo-)K\"ahler form
$\omega_0 = \frac{\im}{2}\delta_{k\cc{l}} \d z^k \wedge
\d\cc{z}^l$. In this case it is well-known that one can define an
associative product on the functions on $\mathbb C^n$ that are
polynomials in the coordinates $(z,\cc{z})$ using the Wick resp.
normal ordering of creation and annihilation operators. Replacing
$\hbar$ by $\frac{\nu}{\im}$ in the resulting formula for the
product of two polynomial functions $F,G$ on $\mathbb C^n$ one
obtains
\BEQ{WickflatEq}
F \starWi G = \sum_{l=0}^\infty
\frac{1}{l!}\left(\frac{2\nu}{\im}\right)^l
\delta^{i_1\cc{j}_1}\ldots
\delta^{i_l\cc{j}_l}
\frac{\partial^l F}{\partial z^{i_1}\cdots\partial z^{i_l}}
\frac{\partial^l G}{\partial \cc{z}^{j_1}\cdots\partial
\cc{z}^{j_l}}.
\EEQ
From this explicit expression it is obvious that one obtains a star
product, the so-called Wick product, also denoted by $\starWi$ on
$(\mathbb C^n,\omega_0)$ using this formula for $f,g\in
\Cinf{\mathbb C^n}[[\nu]]$ instead of $F,G$. Considering the
explicit shape of $\starWi$ one immediately notices that $f \starWi
h$ coincides with the pointwise product of $f$ and $h$ for all
holomorphic functions $h\in\Hol{\mathbb C^n}$ and all
$f\in\Cinf{\mathbb C^n}$. Moreover, we have $h' \starWi g = h' g$
for all anti-holomorphic functions $h'\in \AHol{\mathbb C^n}$ and
all $g\in \Cinf{\mathbb C^n}$.

In order to generalize this property of the star product $\starWi$
to star products on an arbitrary pseudo-K\"ahler manifold the
following definition offers itself: A star product on
$(M,\omega,I)$ is called star product of Wick type if in a local
holomorphic chart the describing bidifferential operators $C_r$ for
$r\geq 1$ have the shape
\BEQ{WickbidiffOplocEq}
C_r (f,g)= \sum_{K,\cc{L}} C_r^{K ;\cc{L}}\frac{\partial^{|K|}
f}{\partial z^K}\frac{\partial^{|\cc{L}|} g}{\partial
\cc{z}^L}
\EEQ
with certain coefficient functions $C_r^{K ;\cc{L}}$. Obviously
this characterization of star products of Wick type is a global
notion on pseudo-K\"ahler manifolds which indeed is equivalent to
the following (cf. \cite[Thm. 4.7]{BorWal96a}):

\begin{DEFINITION}
A star product $\star$ on a pseudo-K\"ahler manifold $(M,\omega,I)$
is called star product of Wick type if for all open subsets
$U\subseteq M$ and all $f,g, h, h'\in \Cinf{M}$ with $h|_U \in
\Hol{U}, h'|_U \in \AHol{U}$ the equations
\BEQ{WickTypAllgDefEq}
f \star h |_U = f h|_U\qquad \textrm{ and } \qquad h' \star g|_U =
h' g|_U
\EEQ
are valid.
\end{DEFINITION}

Compared to the definition of a star product with separation of
variables as introduced by Karabegov the r\^{o}les of holomorphic
and anti-holomorphic directions are interchanged here and according
to our convention for the sign of the Poisson bracket the star
products with separation of variables that are considered by
Karabegov are star products on $(M,-\omega,I)$. Briefly, this means
that a star product $\star$ on $(M,\omega,I)$ is a star product of
Wick type if and only if the opposite star product $\starfl$ which
is defined by $f\starfl g := g \star f$ for $f,g\in
\Cinf{M}[[\nu]]$ is a star product with separation of variables on
$(M,-\omega,I)$. In physics literature the star products with
separation of variables on $(M,\omega,I)$ are usually called star
products of anti-Wick type since in the example $(\mathbb
C^n,\omega_0)$ they arise from the anti-Wick resp. anti-normal
ordering of creation and annihilation operators. In this case an
explicit formula can be obtained very easily and is given by
\BEQ{AntiWickflatEq}
f \starAW g =\sum_{l=0}^\infty
\frac{1}{l!}\left(-\frac{2\nu}{\im}\right)^l
\delta^{i_1\cc{j}_1}\ldots
\delta^{i_l\cc{j}_l}
\frac{\partial^l f}{\partial \cc{z}^{j_1}\cdots\partial
\cc{z}^{j_l}}
\frac{\partial^l g}{\partial z^{i_1}\cdots\partial z^{i_l}}.
\EEQ
This example motivates the definition:
\begin{DEFINITION}
A star product $\star$ on a pseudo-K\"ahler manifold $(M,\omega,I)$
is called star product of anti-Wick type if for all open subsets
$U\subseteq M$ and all $f,g, h, h'\in \Cinf{M}$ with $h|_U \in
\Hol{U}, h'|_U \in \AHol{U}$ the equations
\BEQ{AntiWickTypAllgDefEq}
h \star f |_U = hf |_U\qquad \textrm{ and } \qquad g \star h' |_U =
g h' |_U
\EEQ
are valid.
\end{DEFINITION}
There is a very simple relation between star products of Wick type
and star products of anti-Wick type that can be established using
the parity operator $\Pa$ defined by $\Pa := (-1)^{\deg_\nu}$,
where $\deg_\nu:= \nu \partial_\nu$, that obviously satisfies
$\Pa^2=\id$. For any given star product $\star$ on $(M,\omega,I)$
we consider $\star_{\pa,\opp}$ defined by $f
\star_{\pa,\opp} g := \Pa ((\Pa f)\starfl (\Pa g)) =
\Pa ((\Pa g)\star(\Pa f))$ that again is a star product on
$(M,\omega,I)$. Obviously the mapping $\star \mapsto
\star_{\pa,\opp}$ establishes a bijection on the set of
star products on $(M,\omega,I)$ and from $\Pa^2=\id$ it is evident
that ${\star_{\pa,\opp}}_{\pa,\opp}= \star$. It follows from the
very definitions that $\star$ is of Wick type resp. anti-Wick type
if and only if $\star_{\pa,\opp}$ is of anti-Wick type resp. Wick
type. Due to this relation we shall mostly restrict our attention
to the investigation of star products of Wick type and deduce the
respective analogous statements for the star products of anti-Wick
type using the aforementioned mapping.
\section{Fedosov Star Products, Equivalences and Characteristic
Classes}
\label{FedProdEquCharClaSec}
In this section we shall briefly recall Fedosov's construction of
star products using some different fibrewise products that are
adapted to the special geometric situation on a pseudo-K\"ahler
manifold and are modelled on the examples $\starWi$ resp. $\starAW$
on $\mathbb C^n$. Moreover, we shall relate the obtained star
products to star products that are obtained using the usual
fibrewise Weyl product as in Fedosov's original construction by
explicitly constructing equivalence transformations between the
different products. Using our result in \cite[Thm. 4.4]{Neu99} and
these equivalences it is trivial to compute the characteristic
classes of the constructed star products. In addition we shall
establish some further relations between the obtained star
products.

In the following $(M,\omega,I)$ is a pseudo-K\"ahler manifold and
$\nabla$ denotes the corresponding pseudo-K\"ahler connection. The
basis for Fedosov's construction (cf. \cite{Fed94,Fed96}) is the
following $\mathbb C[[\nu]]$-module
\BEQ{WLDefEq}
\WL(M) := \left(\mathsf{X}_{s=0}^\infty \Ginf{
\mbox{$\bigvee$}^s T^*M \otimes
\mbox{$\bigwedge$} T^*M}\right)[[\nu]].
\EEQ
If there is no possibility for confusion we simply write $\WL$. By
$\WL^k$ we denote the elements of anti-symmetric degree $k$ and set
$\mathcal W:= \WL^0$. For $a,b \in
\WL$ one defines their pointwise (undeformed) product denoted by
$\mu(a\otimes b)= ab$ by the symmetric $\vee$-product in the first
factor and the anti-symmetric $\wedge$-product in the second factor
which turns $\WL$ into a super-commutative associative algebra.
Then the degree-maps $\deg_s$, $\deg_a$, $\deg_\nu:= \nu
\partial_\nu$ and $\Deg:= \deg_s + 2 \deg_\nu$ are defined in the
usual way and yield derivations with respect to $\mu$. For a vector
field $X\in
\Ginf{TM}$ the symmetric resp. anti-symmetric insertion maps are
denoted by $i_s(X)$ resp. $i_a(X)$ which are super-derivations of
symmetric degree $-1$ resp. $0$ and anti-symmetric degree $0$ resp.
$-1$. In local coordinates $x^1, \ldots ,x^{2n}$ for $M$ as a real
$2n$-dimensional manifold we define $\delta:= (1 \otimes
\d x^i)i_s(\partial_i)$ and $\delta^* := (\d x^i \otimes
1)i_a(\partial_i)$ satisfying the relations
$\delta^2=(\delta^*)^2=0$ and $\delta\delta^* + \delta^*\delta =
\deg_s + \deg_a$. Denoting the projection onto the part of
symmetric and anti-symmetric degree $0$ by $\sigma:\WL \to
\Cinf{M}[[\nu]]$ one has
\BEQ{HodZerEq}
\delta \delta^{-1} + \delta^{-1}\delta + \sigma = \id,
\EEQ
where $\delta^{-1} a := \frac{1}{k+l} \delta^* a$ for $\deg_s a = k
a$, $\deg_a a = l a$ with $k+l \neq 0$ and $\delta^{-1}a:=0$ else.

Now we consider different fibrewise associative deformations of the
pointwise product $\mu$: For $a,b \in \WL$ the usual fibrewise Weyl
product $\circF$ which was used in Fedosov's original construction
is defined by
\BEQ{CircFDefEq}
a \circF b:= \mu \circ \exp \left(\frac{\nu}{2} \Lambda^{ij}
i_s(\partial_i) \otimes i_s(\partial_j)\right) (a \otimes b),
\EEQ
where $\Lambda^{ij}$ denotes the components (in real coordinates)
of the Poisson tensor corresponding to the symplectic form $\omega$
that are related to the ones of $\omega$ by the equation
$\omega_{kj}\Lambda^{ij}= \delta^i_k$. Moreover, we can define the
fibrewise Wick product $\circWi$ by
\BEQ{CircWiDefEq}
a \circWi b := \mu \circ \exp \left( \frac{2\nu}{\im}
g^{k\cc{l}}i_s(Z_k)\otimes i_s(\cc{Z}_l)\right)(a\otimes b)
\EEQ
and the fibrewise anti-Wick product by
\BEQ{CircAWDefEq}
a \circAW b :=\mu \circ \exp \left( -\frac{2\nu}{\im}
g^{k\cc{l}}i_s(\cc{Z}_l)\otimes i_s(Z_k) \right)(a\otimes b),
\EEQ
where we have used local holomorphic coordinates of $M$ as a
$n$-dimensional complex manifold and the Poisson tensor is written
as $\Lambda = \frac{2}{\im} g^{k\cc{l}} Z_k \wedge \cc{Z}_l$ with
$g^{k\cc{l}}g_{k\cc{n}}= \delta^{\cc{l}}_{\cc{n}}$. Using these
products we define $\deg_a$-graded super-commutators with respect
to $\circ_{\diamond}$ and set $\ad_{\diamond}(a)b :=
[a,b]_{\circ_{\diamond}}$, where $\diamond$ stands for $\TinyF$,
$\TinyWi$ or $\TinyAW$.

The considered fibrewise products turn out to be fibrewisely
equivalent, i.e. we have
\BEQ{FibEquiEq}
\mathcal S^{-1}\left((\mathcal S a) \circWi (\mathcal S b)\right)=
a\circF b = \mathcal S\left((\mathcal S^{-1} a) \circAW (\mathcal
S^{-1} b)\right) \qquad \forall a,b \in \WL,
\EEQ
where the fibrewise equivalence transformation is given by
\BEQ{SDefEq}
\mathcal S :=\exp \left(\frac{\nu}{\im}\Delta_{\mbox{\rm\tiny fib}}
\right) \qquad \textrm{ with } \qquad \Delta_{\mbox{\rm\tiny fib}}:=
g^{k\cc{l}} i_s(Z_k) i_s(\cc{Z}_l).
\EEQ
Because of the explicit shape of the fibrewise products $\circF$,
$\circWi$, $\circAW$ it is obvious that $\delta$ is a
super-derivation of anti-symmetric degree $1$. In addition one
immediately verifies that the total degree-map $\Deg$ is a
derivation with respect to these products such that $(\WL,
\circ_\diamond)$ is a formally $\Deg$-graded algebra, where again
$\diamond$ stands for $\TinyF$, $\TinyWi$ or $\TinyAW$. Using the
pseudo-K\"ahler connection $\nabla$ on $TM$ that extends in the
usual way to a connection on $T^*M$ and symmeric resp.
anti-symmetric products thereof we define (using the same symbol as
for the connection) the map $\nabla$ on $\WL$ by $\nabla := (1
\otimes \d x^i) \nabla_{\partial_i}$. Since the pseudo-K\"ahler
connection is symplectic $\nabla$ turns out to be a
super-derivation of anti-symmetric degree $1$ and symmetric degree
$0$ with respect to the fibrewise product $\circF$. From the
additional property of the pseudo-K\"ahler connection that $\nabla
I=0$ one directly finds that $[\nabla,\Delta_{\mbox{\rm\tiny
fib}}]=0$ and hence $\nabla =
\mathcal S \nabla \mathcal S^{-1} =
\mathcal S^{-1}\nabla\mathcal S$. These relations together with
equation (\ref{FibEquiEq}) imply that $\nabla$ is also a
super-derivation with respect to the fibrewise products $\circWi$
and $\circAW$. Moreover, $[\delta,\nabla]=0$ since the
pseudo-K\"ahler connection is torsion free and $\nabla^2 =
-\frac{1}{\nu}\adF(R)$, where $R := \frac{1}{4}\omega_{it}
R^t_{jkl} \d x^i \vee \d x^j \otimes \d x^k \wedge \d x^l=
\frac{\im}{2}g_{k \cc{t}}R^{\cc{t}}_{\cc{l}i
\cc{j}} \d z^k \vee \d \cc{z}^l\otimes \d z^i \wedge \d\cc{z}^j \in
\WL^2$ involves the curvature of the connection. As consequences of
the Bianchi identities $R$ satisfies $\delta R =\nabla R=0$. Since
$\nabla$ commutes with $\mathcal S$ we also have $\nabla^2=
\mathcal S \nabla^2 \mathcal S^{-1} = -\frac{1}{\nu}\mathcal S
\adF(R)\mathcal S^{-1} = -\frac{1}{\nu}\adWi(\mathcal S R) =
-\frac{1}{\nu} \adWi(R)$. Here we again have used equation
(\ref{FibEquiEq}) and $\mathcal S R = R +
\frac{\nu}{\im}\Delta_{\mbox{\rm\tiny fib}}R$ which implies the
above equation since $\deg_s\Delta_{\mbox{\rm\tiny fib}}R=0$ and
hence $\adWi(\Delta_{\mbox{\rm\tiny fib}}R)=0$. Analogously one
finds $\nabla^2 = -\frac{1}{\nu} \adAW(R)$.

Now remember the following two theorems which are just slight
generalizations of Fedosov's original theorems in \cite[Thm. 3.2,
3.3]{Fed94}, where $\diamond$ stands for $\TinyF$, $\TinyWi$ or
$\TinyAW$:

\begin{THEOREM}\label{littlerThm}
Let $\Omega_\diamond = \sum_{i=1}^\infty \nu^i
\Omega_{\diamond,i}$ denote a formal series of closed two-forms on
$M$ and let $s_\diamond
=\sum_{k=3}^\infty {s_\diamond}^{(k)} \in \mathcal W$ be given with
$\sigma(s_\diamond)=0$ and $\Deg {s_\diamond}^{(k)} = k
{s_\diamond}^{(k)}$. Then there exists a unique element $r_\diamond
\in \WL^1$ of total degree $\geq 2$ such that
\BEQ{littlerEq}
\delta r_\diamond = \nabla r_\diamond - \frac{1}{\nu}
r_\diamond \circ_\diamond r_\diamond + R + 1\otimes
\Omega_\diamond\quad \textrm{ and }\quad \delta^{-1} r_\diamond =
s_\diamond.
\EEQ
Moreover, $r_\diamond = \sum_{k=2}^\infty {r_\diamond}^{(k)}$ with
$\Deg {r_\diamond}^{(k)} = k {r_\diamond}^{(k)}$ satisfies the
formula
\BEQ{littlerRecEq}
r_\diamond = \delta s_\diamond + \delta^{-1}\left(
\nabla r_\diamond - \frac{1}{\nu}r_\diamond \circ_\diamond
r_\diamond + R + 1\otimes \Omega_\diamond\right)
\EEQ
from which $r_\diamond$ can be determined recursively. In this case
the Fedosov derivation
\BEQ{FedDerDefEq}
\mathfrak D_\diamond := -\delta + \nabla - \frac{1}{\nu}
\ad_\diamond (r_\diamond)
\EEQ
is a super-derivation of anti-symmetric degree $1$ with respect to
$\circ_\diamond$ and has square zero: ${\mathfrak D_\diamond}^2=0$.
\end{THEOREM}

\begin{THEOREM}\label{tauThm}
Let $\mathfrak D_\diamond = -\delta + \nabla -
\frac{1}{\nu}\ad_\diamond(r_\diamond)$ be given as in
(\ref{FedDerDefEq}) with $r_\diamond$ as in (\ref{littlerEq}).
\begin{enumerate}
\item
Then for any $f \in \Cinf{M}[[\nu]]$ there exists a unique element
$\tau_\diamond (f)\in \ker(\mathfrak D_\diamond)\cap \mathcal W$
such that $\sigma(\tau_\diamond (f)) =f$ and $\tau_\diamond :
\Cinf{M}[[\nu]] \to \ker(\mathfrak D_\diamond)\cap \mathcal W$ is
$\mathbb C[[\nu]]$-linear and referred to as the Fedosov-Taylor
series corresponding to $\mathfrak D_\diamond$.
\item
For $f\in \Cinf{M}$ we have $\tau_\diamond (f) =
\sum_{k=0}^\infty\tau_\diamond (f)^{(k)}$, where $\Deg \tau_\diamond
(f)^{(k)}= k\tau_\diamond (f)^{(k)}$ which can be obtained by the
following recursion formula
\BEQ{tauRecEq}
\begin {array} {c}
      \tau_\diamond (f)^{(0)} = f \\
      \displaystyle
      \tau_\diamond (f)^{(k+1)} =
      \delta^{-1} \left(\nabla \tau_\diamond (f)^{(k)} -
      \frac{1}{\nu}
      \sum_{l=0}^{k-1}\ad_\diamond ({r_\diamond}^{(l+2)})
      \tau_\diamond (f)^{(k-l)}\right).
\end {array}
\EEQ
\item
Since $\mathfrak D_\diamond$ is a $\circ_\diamond$-super-derivation
of anti-symmetric degree $1$ as constructed in Theorem
\ref{littlerThm} $\ker(\mathfrak D_\diamond)\cap \mathcal W$ is a
$\circ_\diamond$-subalgebra and a new associative product
$\ast_\diamond$ for $\Cinf{M}[[\nu]]$ is defined by pull-back of
$\circ_\diamond$ via $\tau_\diamond$, which turns out to be a star
product on $(M,\omega,I)$.
\end{enumerate}
\end{THEOREM}

One first result concerning all the star products constructed above
is the following:

\begin{PROPOSITION}\label{VeyTypProp}
The star products $\astF$, $\astWi$ and $\astAW$ on $(M,\omega,I)$
are of Vey type, i.e. the bidifferential operators $C_{\diamond,
r}$ describing the star product in order $r$ of the formal
parameter via $f \ast_\diamond g = \sum_{r=0}^\infty \nu^r
C_{\diamond,r}(f,g)$ for $f,g\in \Cinf{M}$ are of order $(r,r)$.
\end{PROPOSITION}
\begin{PROOF}
We write the term of total degree $k\neq 0$ in the Fedosov-Taylor
series as $\tau_\diamond(f)^{(k)}
= \sum_{l=0}^{[\frac{k-1}{2}]}
\nu^l{\tau_\diamond(f)}_{k-2l}^{(k)}$ where
${\tau_\diamond(f)}_{k-2l}^{(k)} \in
\Ginf{\bigvee^{k-2l}T^*M}$ for $f \in \Cinf{M}$. A straightforward
computation yields that
\[
f \ast_\diamond g = fg + \sum_{r=1}^\infty \nu^r
\sum_{l=0}^{r-1}\sum_{k=0}^l \mu_{\diamond,r-l}(
{\tau_\diamond(f)}_{r-l}^{(r-l+2k)},
{\tau_\diamond(g)}_{r-l}^{(r+l-2k)} ),
\]
where we have written $\mu_{\diamond,k}({}\cdot{},{}\cdot{})$ for
the fibrewise bidifferential operator occurring at order $k\geq1$
of the formal parameter in the fibrewise product $\circ_\diamond$.
Now a lengthy but straightforward induction on the total degree
yields:
\begin{LEMMA}
For all $k\geq 1$ and $0\leq l \leq[\frac{k-1}{2}]$ the mapping
\BEQ{tauDiffOpEq}
\Cinf{M}\ni f \mapsto {\tau_\diamond(f)}_{k-2l}^{(k)} \in
\Ginf{\mbox{$\bigvee$}^{k-2l}T^*M}
\EEQ
is a differential operator of order $k-l$.
\end{LEMMA}
According to this lemma ${\tau_\diamond}_{r-l}^{(r-l+2k)}$ is a
differential operator of order $r-l+k$ and
${\tau_\diamond}_{r-l}^{(r+l-2k)}$ is a differential operator of
order $r-k$. Observing the ranges of summation in the expression
for $C_{\diamond,r}$ given above it is easy to see that the highest
occurring order of differentiation is $r$ in both arguments,
proving the assertion.
\end{PROOF}

In fact an analogous statement to the one of the above proposition
which shall be published elsewhere \cite{Neu02b} holds in even more
general situations, where one can construct Fedosov-like products.

Another observation that is indeed rather technical but shall turn
out to be useful and important in the next section is that one can
restrict to elements $s_\diamond$ that contain no part of symmetric
degree $1$ without omitting any of the possible star products in
case one allows for all possible formal series of closed two-forms
$\Omega_\diamond$.

\begin{LEMMA}\label{NormCondLem}
Let $\mathfrak D_\diamond$ and $\mathfrak D'_\diamond$ be two
Fedosov derivations for $\circ_\diamond$ constructed from different
data $(\Omega_\diamond, s_\diamond)$ and
$(\Omega'_\diamond,s'_\diamond)$, then we have
\BEQ{DiffDatEquiEq}
\mathfrak D_\diamond = \mathfrak D'_\diamond \iff
r'_\diamond - r_\diamond = 1 \otimes B_\diamond \iff
\begin{array}{l} s'_\diamond - s_\diamond = B_\diamond \otimes 1
\quad
\textrm {and}\\
\Omega_\diamond - \Omega'_\diamond = \d B_\diamond,
\end{array}
\EEQ
where $B_\diamond = \sum_{i=1}^\infty \nu^i B_{\diamond,i} \in
\Ginf{T^*M}[[\nu]]$ is a formal series of one-forms on $M$. In case
one of the conditions in (\ref{DiffDatEquiEq}) is fulfilled
obviously the star products $\ast_\diamond$ and $\ast'_\diamond$
coincide.
\end{LEMMA}
\begin{PROOF}
Obviously $\mathfrak D_\diamond = \mathfrak D'_\diamond \iff
\frac{1}{\nu}\ad_\diamond(r'_\diamond - r_\diamond)=0$ but
this is equivalent to $r'_\diamond - r_\diamond$ being central and
hence $r'_\diamond - r_\diamond = 1 \otimes B_\diamond$ with a
formal series $B_\diamond$ of one-forms that vanishes at order zero
of the formal parameter since the total degree of $r'_\diamond$ and
$r_\diamond$ is at least $2$. For the proof of the second
equivalence first let $r'_\diamond - r_\diamond = 1 \otimes
B_\diamond$ from which we obtain applying $\delta^{-1}$ to this
equation and using (\ref{littlerEq}) that $s'_\diamond - s_\diamond
= B_\diamond \otimes 1$. To obtain the second equation we apply
$\delta$ to $r'_\diamond - r_\diamond = 1 \otimes B_\diamond$ and a
straightforward computation using (\ref{littlerEq}) and
$\nabla(1\otimes B_\diamond)= 1 \otimes \d B_\diamond$ yields $0 =
\d B_\diamond+\Omega'_\diamond-\Omega_\diamond$ proving the first
direction of the second equivalence. Now assume that the relations
between the data $(\Omega_\diamond, s_\diamond)$ and
$(\Omega'_\diamond, s'_\diamond)$ are satisfied and define
$A_\diamond:=r_\diamond - r'_\diamond + 1 \otimes B_\diamond$ then
an easy computation using (\ref{littlerEq}) yields that
$A_\diamond$ fulfills the equations
\[
-\delta A_\diamond + \nabla A_\diamond - \frac{1}{\nu}
\ad_\diamond (r_\diamond) A_\diamond +
\frac{1}{\nu} A_\diamond \circ_\diamond A_\diamond =0\quad
\textrm{ and } \quad \delta^{-1}A_\diamond =0.
\]
From these equations one obtains with $\sigma(A_\diamond)=0$ using
the decomposition $\delta\delta^{-1}+\delta^{-1}\delta + \sigma
=\id$ that $A_\diamond$ is a fixed point of the mapping $L$ defined
by $La := \delta^{-1}(\nabla a - \frac{1}{\nu}\ad_\diamond
(r_\diamond)a + \frac{1}{\nu}a \circ_\diamond a)$ for $a\in \WL^1$
which is at least of total degree $2$. Now $L$ raises the total
degree at least by $1$ and using Banach's fixed point theorem (cf.
Corollary \ref{FedConAbbCor} ii.)) we obtain that $L$ has a unique
fixed point. But as well $A_\diamond$ as (trivially) $0$ are fixed
points of $L$ yielding $A_\diamond =0$ by uniqueness and hence
$r'_\diamond - r_\diamond = 1 \otimes B_\diamond$.
\end{PROOF}

The statement of the preceding lemma now means the following:
Starting with a formal series $\Omega_\diamond$ of closed two-forms
and an element $s_\diamond$ with $\sigma(s_\diamond)=0$ that in
addition contains no part of symmetric degree $1$ and passing over
to the normalization condition $\widetilde{s}_\diamond=s_\diamond +
B_\diamond\otimes 1$ is equivalent to passing over to
$\Omega_\diamond + \d B_\diamond$ sticking to the normalization
condition $s_\diamond$. Hence we may without loss of generality
restrict to normalization conditions shaped like $s_\diamond$
provided that we allow for arbitrarily varying closed two-forms
$\Omega_\diamond$.

Now we shall compute Deligne's characteristic class of the star
products $\astF$, $\astWi$ and $\astAW$ constructed above showing
that one can obtain star products of every given equivalence class
with each of the above fibrewise products $\circF$, $\circWi$ and
$\circAW$ by suitable choice of the closed two-forms $\OF$, $\OWi$
and $\OAW$. For details on the definition of Deligne's
characteristic class and methods to compute it the reader is
referred to \cite{Del95,GR99,Gut00,Neu99}.

\begin{PROPOSITION}\label{CharKlaProp}
Consider the star products constructed according to Theorem
\ref{tauThm} iii.).
\begin{enumerate}
\item
Deligne's characteristic classes $c(\astF), c(\astWi),c(\astAW)\in
\frac{[\omega]}{\nu}+ H^2_{\mbox{\rm\tiny dR}}(M,\mathbb C)[[\nu]]$
of the star products $\astF,\astWi,\astAW$ on $(M,\omega,I)$ are
given by
\BEQ{classEq}
c(\astF) = \frac{[\omega]}{\nu} + \frac{[\OF]}{\nu},\quad
c(\astWi)=\frac{[\omega]}{\nu} +
\frac{[\OWi]}{\nu}-\frac{[\rho]}{\im},\quad
c(\astAW)= \frac{[\omega]}{\nu} + \frac{[\OAW]}{\nu} +
\frac{[\rho]}{\im},
\EEQ
where $\rho$ denotes the Ricci form which is the closed two-form
given by $\rho = - \frac{\im}{2} R_{i\cc{j}}\d z^i \wedge
\d\cc{z}^j$ resp. $\Delta_{\mbox{\rm\tiny fib}} R = 1 \otimes \rho$
and $R_{i\cc{j}}$ denotes the components of the Ricci tensor in a
local holomorphic chart.
\item
Because of the properties of the characteristic class i.) implies
that two star products $\ast_\diamond$ and $\ast'_\diamond$
constructed from different data $(\Omega_\diamond,s_\diamond)$ and
$(\Omega'_\diamond,s'_\diamond)$ using the same fibrewise product
$\circ_\diamond$ are equivalent if and only if
$[\Omega_\diamond]=[\Omega'_\diamond]$. In this case
$\Omega_\diamond-\Omega'_\diamond = \d C_\diamond$ with a formal
series $C_\diamond = \sum_{i=1}^\infty \nu^i C_{\diamond,i}\in
\Ginf{T^*M}[[\nu]]$ of one-forms on $M$ and an equivalence
transformation can be constructed in the following way: There is a
uniquely determined element $h_\diamond = \sum_{k=3}^\infty
h_\diamond^{(k)}\in\mathcal W$ with $\Deg h_\diamond^{(k)} = k
h_\diamond^{(k)}$ and $\sigma(h_\diamond) =0$ such that
\BEQ{littlehEq}
r'_\diamond - r_\diamond -
\frac{\exp\left(\frac{1}{\nu}\ad_\diamond(h_\diamond)\right)-\id}
{\frac{1}{\nu}\ad_\diamond(h_\diamond)}\left(\mathfrak D_\diamond
h_\diamond\right) = 1 \otimes C_\diamond.
\EEQ
For this $h_\diamond$ one has $\mathfrak D'_\diamond =
\exp\left(\frac{1}{\nu}\ad_\diamond(h_\diamond)\right)
\mathfrak D_\diamond \exp \left(-
\frac{1}{\nu}\ad_\diamond(h_\diamond)\right)$ and
$A_{h_{\diamond}}$ defined by
\BEQ{AhaDefEq}
A_{h_{\diamond}}f :=
\sigma\left(\exp\left(\frac{1}{\nu}\ad_\diamond(h_\diamond)\right)
\tau_\diamond(f)\right)
\EEQ
for $f \in \Cinf{M}[[\nu]]$ the inverse of which is given by
$A^{-1}_{h_{\diamond}}f = \sigma
(\exp\left(-\frac{1}{\nu}\ad_\diamond(h_\diamond)\right)
\tau'_\diamond(f))$
is an equivalence from $\ast_\diamond$ to $\ast'_\diamond$, i.e.
$A_{h_{\diamond}}(f
\ast_\diamond g) = (A_{h_{\diamond}} f)\ast'_\diamond
(A_{h_{\diamond}}g) \,\forall f,g\in \Cinf{M}[[\nu]]$.
\end{enumerate}
\end{PROPOSITION}
\begin{PROOF}
The characteristic class of $\astF$ has been determined in
\cite[Thm. 4.4]{Neu99} and is given by the above expression. We
consider a Fedosov derivation $\DWi$ for $\circWi$ then obviously
$\mathfrak D'_\TinyF:=\mathcal S^{-1} \DWi \mathcal S = -\delta +
\nabla- \frac{1}{\nu}\adF(\mathcal S^{-1}\rWi)$ defines a Fedosov
derivation for $\circF$. Now we define $s'_\TinyF:=
\delta^{-1}\mathcal S^{-1}\rWi$ and observe that $\sigma(s'_\TinyF)
=0$ and that $s'_\TinyF$ contains no terms of total degree lower
than $3$. Moreover, the equation that is satisfied by $\rWi$ yields
that $r'_\TinyF:=\mathcal S^{-1}\rWi$ obeys the equation $\delta
r'_\TinyF = \nabla r'_\TinyF - \frac{1}{\nu}r'_\TinyF
\circF r'_\TinyF + R
-\frac{\nu}{\im} 1 \otimes \rho + 1 \otimes \OWi$ and hence from
the first statement of the proposition we get $c(\ast'_\TinyF) =
\frac{[\omega]}{\nu} + \frac{[\OWi]}{\nu}-\frac{[\rho]}{\im}$. But
this implies the result on the characteristic class of $\astWi$ as
it is straightforward using the relation between $\DWi$ and
$\mathfrak D'_\TinyF$ to see that $S$ defined by $Sf := \sigma
(\mathcal S \tau'_\TinyF(f))$ for $f\in \Cinf{M}[[\nu]]$ yields an
equivalence transformation from $\ast'_\TinyF$ to $\astWi$ (cf.
\cite[Prop. 2.3]{BorNeu98a}) implying that their characteristic
classes coincide. For the proof of the formula for $c(\astAW)$ one
proceeds completely analogously replacing $\mathcal S$ by $\mathcal
S^{-1}$ in the above argument causing the opposite sign in front of
the Ricci form $\rho$. For the proof of ii.) first note that two
star products are equivalent if and only if their characteristic
classes coincide yielding that $\ast_\diamond$ and $\ast'_\diamond$
are equivalent if and only if
$[\Omega_\diamond]=[\Omega'_\diamond]$ by i.). For the explicit
construction of the equivalence transformation similar to the one
in \cite[Thm. 5.2]{BorNeu98a} one finds by easy computation that
$\mathfrak D'_\diamond =
\exp\left(\frac{1}{\nu}\ad_\diamond(h_\diamond)\right)
\mathfrak D_\diamond \exp \left(-
\frac{1}{\nu}\ad_\diamond(h_\diamond)\right)$ if and only if
(\ref{littlehEq}) is satisfied. Applying $\mathfrak D_\diamond$ to
this equation and using (\ref{littlerEq}) together with $\mathfrak
D_\diamond^2=0$ one finds that the necessary condition for
(\ref{littlehEq}) to be solvable is
$\Omega_\diamond-\Omega'_\diamond = \d C_\diamond$ which is
fulfilled by assumption. Solving (\ref{littlehEq}) for $\mathfrak
D_\diamond h_\diamond$ and applying $\delta^{-1}$ to the obtained
equation one gets using $\delta\delta^{-1}+ \delta^{-1}\delta +
\sigma = \id$ and $\sigma (h_\diamond)=0$ that $h_\diamond$ has to
satisfy the equation
\BEQ{hdiamRecEq}
h_\diamond = C_\diamond \otimes 1 + \delta^{-1}\left(
\nabla h_\diamond - \frac{1}{\nu}\ad_\diamond (r_\diamond)
h_\diamond - \frac{\frac{1}{\nu}\ad_\diamond(h_\diamond)}{
\exp(\frac{1}{\nu}\ad_\diamond(h_\diamond))-\id}(r'_\diamond-
r_\diamond)
\right).
\EEQ
Taking into account the total degrees of the involved mappings it
is easy to see that the right-hand side of this equation is of the
form $L h_\diamond$, where $L$ is a contracting mapping and again
from Banach's fixed point theorem (cf. Corollary \ref{FedConAbbCor}
i.)) we get the existence and uniqueness of the element
$h_\diamond$ satisfying $h_\diamond = L h_\diamond$. Now it remains
to prove that this fixed point $h_\diamond$ actually solves
(\ref{littlehEq}). To this end one defines
$A_\diamond:=\frac{\frac{1}{\nu}\ad_\diamond(h_\diamond)}{
\exp(\frac{1}{\nu}\ad_\diamond(h_\diamond))-\id}(r'_\diamond-
r_\diamond)- \mathfrak D_\diamond h_\diamond - 1 \otimes
C_\diamond$ and applying $\mathfrak D_\diamond$ to this equation
one derives an equation for $A_\diamond$ of the shape $\mathfrak
D_\diamond A_\diamond =
R_{h_\diamond,r'_\diamond,r_\diamond}(A_\diamond)$, where
$R_{h_\diamond,r'_\diamond,r_\diamond}$ is a linear mapping that
does not decrease the total degree. The explicit shape of
$R_{h_\diamond,r'_\diamond,r_\diamond}$ is of minor importance.
Using $\delta^{-1}A_\diamond
=\sigma(A_\diamond)=0$ this yields that $A_\diamond$ satisfies
\[
A_\diamond = \delta^{-1}\left(\nabla A_\diamond - \frac{1}{\nu}
\ad_\diamond (r_\diamond)A_\diamond - R_{h_\diamond,r'_\diamond,
r_\diamond}(A_\diamond)\right).
\]
Again this is a fixed point equation with a unique solution, but as
$R_{h_\diamond,r'_\diamond,r_\diamond}$ is linear $0$ also solves
this fixed point equation implying $A_\diamond=0$ and hence
$h_\diamond$ with $L h_\diamond = h_\diamond$ actually solves
(\ref{littlehEq}). Now as in i.) it is straighforward to conclude
from $\mathfrak D'_\diamond = \exp\left(\frac{1}{\nu}
\ad_\diamond(h_\diamond)\right) \mathfrak D_\diamond \exp
\left(-\frac{1}{\nu}\ad_\diamond(h_\diamond)\right)$ that
$A_{h_{\diamond}}$ defined in (\ref{AhaDefEq}) is an equivalence
transformation from $\ast_\diamond$ to $\ast'_\diamond$.
\end{PROOF}

Using the preceding proposition we are in the position to decide
under which conditions on the respective data $(\sF,\OF)$,
$(\sWi,\OWi)$ and $(\sAW,\OAW)$ the star products $\astF$, $\astWi$
and $\astAW$ are equivalent and in this case we can moreover
explicitly give an equivalence transformation using suitable
combinations of the constructed transformations given in the proof
of the proposition.

To conclude this section we shall now give a condition on the data
$(\sWi,\OWi)$ and $(\sAW,\OAW)$ under which we have that $\astAW =
{(\astWi)}_{\pa,\opp}$ which shall enable us to restrict all our
further considerations to star products of Wick type, since from
this relation we can immediately deduce the corresponding statement
about the corresponding star product of anti-Wick type.

\begin{LEMMA}\label{AWWiRelLem}
For $\sAW = \Pa \sWi$ and $\OAW = \Pa \OWi$ we have $\astAW
={(\astWi)}_{\pa,\opp}$, i.e. for all $f,g \in \Cinf{M}[[\nu]]$ we
have
\BEQ{astAWastWiPaOppEq}
f \astAW g = \Pa ((\Pa g)\astWi (\Pa f)).
\EEQ
\end{LEMMA}
\begin{PROOF}
First we observe that $\Pa(a \circWi b ) = (-1)^{|a||b|}(\Pa
b)\circAW (\Pa a)$ because of the explicit shape of the fibrewise
products. Using equation (\ref{littlerEq}) for $\rWi$ one obtains
that $\delta^{-1}\Pa \rWi = \Pa \sWi$ and $\delta\Pa \rWi = \nabla
\Pa \rWi - \frac{1}{\nu}(\Pa \rWi) \circAW (\Pa \rWi) + R + 1
\otimes \Pa \OWi$ and hence for $\OAW = \Pa \OWi$ and $\sAW = \Pa
\sWi$ we have $\rAW = \Pa \rWi$ implying that $\DAW \Pa = \Pa
\DWi$. From this relation it is easy to deduce that $\Pa \tauWi(f)
=\tauAW (\Pa f)$ for all $f\in\Cinf{M}[[\nu]]$ and using this
equation the proof of (\ref{astAWastWiPaOppEq}) is straightforward.
\end{PROOF}

\section{Unique Characterization of the Star Products of Wick Type}
\label{WickCharSec}
In this section we shall uniquely characterize the star products of
Wick type among the star products $\astWi$ constructed in the
preceding section, which actually not all are of Wick type although
the fibrewise product $\circWi$ used to obtain them is of this
shape. Throughout the whole section we shall assume that $\DWi$ is
a Fedosov derivation constructed as in Theorem \ref{littlerThm}
using an element $\sWi$ as normalization condition for $\rWi$ that
contains no part of symmetric degree $1$, which according to Lemma
\ref{NormCondLem} represents no restriction for the obtained star
products. Throughout this section we shall extensively make use of
the notations and results given in Appendix \ref{ComStrConSec} that
are crucial for our investigations.

In \cite{BorWal96a} Bordemann and Waldmann have considered the star
product $\astWi$ in case $\sWi=\OWi=0$ and have shown that the
resulting star product is of Wick type. In their proof some special
property of the element $\rWi$ determining the Fedosov derivation
was crucial. In this special case one has according to \cite[Lemma
4.5]{BorWal96a} for all $q,p\geq 0$ that
\BEQ{rBorWalEq}
\pi_s^{p,0} \rWi = \pi_s^{0,q} \rWi =0.
\EEQ
Particularly this means that in this case $\rWi$ contains no part
of symmetric degree $1$. But in case $\OWi\neq 0$ one immediately
finds that $\rWi^{(2k+1)}$ contains for some $k\geq 1$ a summand of
the shape $\nu^k\delta^{-1}(1\otimes
\Omega_{\TinyWi,k})$. On the other hand one cannot do without
allowing for $\OWi\neq 0$ as to obtain star products of arbitrary
characteristic class because of Proposition \ref{CharKlaProp}. The
only way out can consist in showing that a weaker condition on
$\rWi$ than (\ref{rBorWalEq}) which can also be satisfied for
$\OWi\neq 0$ suffices to guarantee that the resulting star product
is of Wick type. The results of Karabegov on the classification of
star products of Wick type now suggest that the case, where $\OWi$
is of type $(1,1)$ should be of special interest.

\begin{LEMMA}\label{rWickProjLem}
Let $\rWi\in \WL^1$ be constructed as in Theorem \ref{littlerThm}
from $\sWi$ containing no part of symmetric degree $1$, then we
have:
\BEQAR{pizrWiNullEq}
\pi_z \rWi = 0 & \iff & \pi_z(1 \otimes \OWi )=0 \quad\textrm{and}
\quad\pi_z\sWi=0,\\ \label{picczrWiNullEq}
\pi_{\cc{z}} \rWi = 0 & \iff & \pi_{\cc{z}}(1 \otimes \OWi )=0
\quad\textrm{and} \quad\pi_{\cc{z}}\sWi=0.
\EEQAR
Hence $\pi_z\rWi= \pi_{\cc{z}}\rWi=0$ if and only if $\OWi$ is of
type $(1,1)$ and $\pi_z \sWi=\pi_{\cc{z}}\sWi =0$.
\end{LEMMA}
\begin{PROOF}
We first show the equivalence (\ref{pizrWiNullEq}) and first let
$\pi_z(1\otimes \OWi)=\pi_z \sWi =0$. Applying $\pi_z$ to the
equations that determine $\rWi$ we obtain using the formulas given
in Appendix \ref{ComStrConSec} that $\pi_z \rWi$ satisfies the
equations
\[
\delta_z \pi_z \rWi = \nabla_z \pi_z \rWi - \frac{1}{\nu}\pi_z(
(\pi_z\rWi)\circWi \rWi)\quad \textrm{ and }\quad
{\delta_z}^{\!\!\!-1}\pi_z \rWi=0.
\]
Using ${\delta_z}^{\!\!\!-1} \delta_z + \delta_z
{\delta_z}^{\!\!\!-1} + \pi_{\cc{z}}=\id$ these equations yield
after some easy manipulations that $\pi_z \rWi$ is a fixed point of
the mapping $L : \pi_z(\WL^1)\ni a
\mapsto {\delta_z}^{\!\!\!-1}(\nabla_z a - \frac{1}{\nu}
\pi_z(\adWi(\rWi)a))\in \pi_z(\WL^1)$ that raises the
total degree at least by $1$ and hence has a unique fixed point
according to Corollary \ref{FedConAbbCor} iii.). Obviously $L$ has
$0$ as trivial fixed point implying that $\pi_z
\rWi =0$ by uniqueness. For the other direction of the equivalence
let $\pi_z \rWi=0$ and again apply $\pi_z$ to the equations that
determine $\rWi$ then it is very easy to see that $\pi_z \rWi =0$
implies that $\pi_z (1\otimes \OWi )=0$ and $\pi_z \sWi=0$. For the
proof of (\ref{picczrWiNullEq}) one proceeds completely
analogously.
\end{PROOF}

An important consequence of the presence of the complex structure
on $M$ is that the explicit expression for $f\astWi g$ does not
depend on the complete Fedosov-Taylor series of $f$ and $g$ but
only on the holomorphic part of $\tauWi(f)$ and the
anti-holomorphic part of $\tauWi(g)$. To see this fact we just
observe that the properties of $\pi_z$ and $\pi_{\cc{z}}$ given in
Appendix \ref{ComStrConSec} imply that
\BEQ{astpizpicczEq}
f \astWi g = \sigma(\tauWi(f)\circWi \tauWi(g))=
\sigma((\pi_z\tauWi(f))\circWi (\pi_{\cc{z}}\tauWi(g))).
\EEQ
This observation raises the question whether it is possible to find
more simple recursion formulas for $\pi_z \tauWi(f)$ and
$\pi_{\cc{z}}\tauWi(g)$ than the one for the whole Fedosov-Taylor
series, since these would completely determine the star product
$\astWi$. Under the preconditions $\pi_z \rWi =0$ and
$\pi_{\cc{z}}\rWi=0$ this question can be answered positively as we
state it in the following proposition.

\begin{PROPOSITION}\label{tauWickeinfachProp}
Let $\rWi\in \WL^1$ be constructed as in Theorem \ref{littlerThm}
from $\sWi$ containing no part of symmetric degree $1$ and let
$\tauWi$ denote the corresponding Fedosov-Taylor series according
to Theorem \ref{tauThm}.
\begin{enumerate}
\item
If $\pi_z \rWi =0$ then $\pi_z\tauWi(f)$ satisfies the equation
\BEQ{piztauWiEq}
\delta_z \pi_z\tauWi(f) = \nabla_z \pi_z \tauWi(f) + \frac{1}{\nu}
\pi_z((\pi_z\tauWi(f))\circWi \rWi)
\EEQ
for $f\in\Cinf{M}[[\nu]]$. In this case $\pi_z \tauWi(f)$ is
uniquely determined from this equation since $\sigma
(\pi_z\tauWi(f))=f$ and can be computed recursively for $f\in
\Cinf{M}$ from
\BEQ{piztauRekEq}
\pi_z \tauWi(f) = f + {\delta_z}^{\!\!\!-1}\left(\nabla_z \pi_z\tauWi(f)
+ \frac{1}{\nu}\pi_z((\pi_z\tauWi(f))\circWi \rWi)\right).
\EEQ
\item
If $\pi_{\cc{z}} \rWi =0$ then $\pi_{\cc{z}}\tauWi(f)$ satisfies
the equation
\BEQ{piccztauWiEq}
\delta_{\cc{z}} \pi_{\cc{z}}\tauWi(f) = \nabla_{\cc{z}}
\pi_{\cc{z}} \tauWi(f) - \frac{1}{\nu}
\pi_{\cc{z}}(\rWi \circWi (\pi_{\cc{z}}\tauWi(f)))
\EEQ
for $f\in\Cinf{M}[[\nu]]$. In this case $\pi_{\cc{z}}
\tauWi(f)$ is uniquely determined from this equation since $\sigma
(\pi_{\cc{z}}\tauWi(f))=f$ and can be computed recursively for
$f\in \Cinf{M}$ from
\BEQ{piccztauRekEq}
\pi_{\cc{z}} \tauWi(f) = f + {\delta_{\cc{z}}}^{\!\!\!-1}
\left(\nabla_{\cc{z}} \pi_{\cc{z}}\tauWi(f)
- \frac{1}{\nu}\pi_{\cc{z}}(\rWi \circWi(\pi_{\cc{z}}\tauWi(f)))
\right).
\EEQ
\end{enumerate}
\end{PROPOSITION}
\begin{PROOF}
Applying $\pi_z$ to the equation $\DWi \tauWi(f)=0$ that is valid
for all $f\in \Cinf{M}[[\nu]]$ we obtain using $\pi_z \rWi=0$ that
\[
\delta_z \pi_z \tauWi(f) = \nabla_z \pi_z\tauWi(f) + \frac{1}{\nu}
\pi_z ((\pi_z \tauWi(f))\circWi \rWi)
\]
since $\pi_z (\rWi \circWi
\tauWi(f))=\pi_z((\pi_z\rWi)\circWi\tauWi(f))=0$. To this equation
we apply ${\delta_z}^{\!\!\!-1}$ and with equation
(\ref{HodgeholaholZerEq}) this yields the stated equation
(\ref{piztauRekEq}). Now the mapping $T_f :\pi_z (\mathcal W)\to
\pi_z(\mathcal W)$ defined by $T_f a := f +{\delta_z}^{\!\!\!-1}(
\nabla_z a + \frac{1}{\nu}\pi_z(a\circWi \rWi))$ turns out to be
contracting and hence according to Corollary \ref{FedConAbbCor}
iii.) there is a unique fixed point $b_f\in
\pi_z(\mathcal W)$ with $T_f b_f = b_f$ that obviously satisfies
$\sigma(b_f)=f$. It remains to show that this fixed point actually
solves the equation $\delta_z b_f = \nabla_z b_f +
\frac{1}{\nu}\pi_z(b_f \circWi \rWi)$ since then $b_f =\pi_z
\tauWi(f)$. To this end one defines $A:= -\delta_z b_f +
\nabla_z b_f + \frac{1}{\nu}\pi_z(b_f \circWi \rWi)$ and an easy
computation using Lemma \ref{SquareLem} yields that $A$ satisfies
$\delta_z A = \nabla_z A - \frac{1}{\nu} \pi_z(A \circWi \rWi)$
from which we obtain again using (\ref{HodgeholaholZerEq}) that $A$
is a fixed point of the mapping $L: \pi_z(\WL^1)\ni a\mapsto
{\delta_z}^{\!\!\!-1}(\nabla_z a - \frac{1}{\nu}\pi_z(a\circWi
\rWi)) \in \pi_z (\WL^1)$. But again this fixed point is unique
and so we get $A=0$ implying that $b_f$ solves the equation that is
solved by $\pi_z\tauWi(f)$ proving that $\pi_z\tauWi(f)$ is
uniquely determined by (\ref{piztauWiEq}) and can be obtained
recursively from (\ref{piztauRekEq}). For the proof of ii.) one
proceeds completely analogously.
\end{PROOF}

An easy consequence of the preceding proposition is:
\begin{LEMMA}\label{tauHolAholLem}
Let $U \subseteq M$ be an open subset of $M$.
\begin{enumerate}
\item
If $\pi_z \rWi=0$ then we have for all $h' \in \Cinf{M}$ with
$h'|_U\in \AHol{U}$ that
\BEQ{tauaholEq}
\pi_z\tauWi(h')|_U = h'|_U.
\EEQ
\item
If $\pi_{\cc{z}}\rWi =0$ then we have for all $h \in \Cinf{M}$ with
$h|_U \in \Hol{U}$ that
\BEQ{tauholeq}
\pi_{\cc{z}}\tauWi(h)|_U = h|_U.
\EEQ
\end{enumerate}
\end{LEMMA}
\begin{PROOF}
To prove (\ref{tauaholEq}) it suffices to show that $h'|_U$
satisfies the equation that determines $\pi_z\tauWi(h')|_U$. But
this is an easy computation using that $\nabla_z h'|_U =0$ since
$h'\in \AHol{U}$ and again that $\pi_z \rWi =0$. The proof of ii.)
is completely analogous.
\end{PROOF}

With the results collected up to now we are already in the position
to give a sufficient condition on which the star product $\astWi$
is of Wick type.

\begin{PROPOSITION}\label{SuffForWickProp}
Let $\DWi$ be the Fedosov derivation constructed from the data
$(\OWi,\sWi)$ and let $\astWi$ be the resulting star product on a
pseudo-K\"ahler manifold $(M,\omega,I)$. Moreover, let $U\subseteq
M$ be an open subset of $M$.
\begin{enumerate}
\item
If $\pi_z(1\otimes
\OWi )=\pi_z(\sWi)=0$ then we have for all $g,h' \in\Cinf{M}$
with $h'|_U\in \AHol{U}$ that
\BEQ{astWiaholEq}
h' \astWi g|_U = h' g|_U.
\EEQ
\item
If $\pi_{\cc{z}}(1\otimes
\OWi )=\pi_{\cc{z}}(\sWi)=0$ then we have for all $f,h
\in\Cinf{M}$ with $h|_U\in \Hol{U}$ that
\BEQ{astWiholEq}
f \astWi h|_U = f h|_U.
\EEQ
\end{enumerate}
Consequently the star product $\astWi$ is of Wick type in case we
have $\pi_z\sWi=\pi_{\cc{z}}\sWi=0$ and $\OWi$ is of type $(1,1)$.
\end{PROPOSITION}
\begin{PROOF}
i.) According to Lemma \ref{rWickProjLem} the conditions on $\sWi$
and $\OWi$ imply that $\pi_z \rWi =0$. But with Lemma
\ref{tauHolAholLem} i.) this directly implies the equation
(\ref{astWiaholEq}). The proof of ii.) is analogous.
\end{PROOF}

\begin{REMARK}
It is remarkable that the statement of the above proposition shows
that the two characterizing properties of a star product of Wick
type can in fact be realized separately. In case equation
(\ref{astWiaholEq}) is satisfied it is easy to see that the
bidifferential operators describing $\astWi$ in a local holomorphic
chart are of the shape $C_r (f,g)=\sum_{K,\cc{L},M}
C_r^{K;\cc{L},M}\frac{\partial^{|K|} f}{\partial
z^K}\frac{\partial^{|\cc{L}|+|M|} g}{\partial
\cc{z}^L\partial z^M}$ for $r\geq 1$. Analogously (\ref{astWiholEq})
implies that we have $C_r(f,g)= \sum_{K,\cc{M},\cc{L}}
C_r^{K,\cc{M};\cc{L}}
\frac{\partial^{|K|+|\cc{M}|} f}{\partial
z^K\partial \cc{z}^M}\frac{\partial^{|\cc{L}|} g}{\partial
\cc{z}^L}$.
\end{REMARK}

Our next aim is to see that the conditions given in the above
proposition are not only sufficient but even necessary for the
equations (\ref{astWiaholEq}) and (\ref{astWiholEq}) to hold. To
prove this we need the following two lemmas.

\begin{LEMMA}
Let $\tauWi$ denote the Fedosov-Taylor series corresponding to
$\DWi$ and let $U\subseteq M$ be an open subset of $M$.
\begin{enumerate}
\item
If we have for all $h'\in\Cinf{M}$ with $h'|_U \in \AHol{U}$ that
the Fedosov-Taylor series satisfies $\pi_z \tauWi(h')|_U =h'|_U$
then we have $\pi_z \rWi =0$.
\item
If we have for all $h\in\Cinf{M}$ with $h|_U \in \Hol{U}$ that the
Fedosov-Taylor series satisfies $\pi_{\cc{z}} \tauWi(h)|_U =h|_U$
then we have $\pi_{\cc{z}} \rWi =0$.
\end{enumerate}
\end{LEMMA}
\begin{PROOF}
We apply $\pi_z$ to the equation $\DWi \tauWi(h')=0$ and restrict
the obtained equation to $U$ and get using $\pi_z\tauWi(h')|_U =
h'|_U$ and considering the term of total degree $l\geq 1$ that
\[
0 = \frac{1}{\nu}\sum_{r=1}^l\pi_z((\pi_z \rWi^{(l+2-r)})\circWi
\tauWi(h')^{(r)})|_U
\]
for all $l\geq 1$. Now we shall show by induction on the total
degree that this equation actually implies $\pi_z \rWi=0$. Let
$(z,V)$ be a local holomorphic chart of $M$ and let $\chi\in
\Cinf{M}$ be a function with $\mathrm{supp}(\chi)\subseteq V$ and
$\chi|_U=1$, where $U\subseteq V$ then for all $1 \leq k \leq
\dim_{\mathbb C}(M)$ the function $\chi \cc{z}^k \in \Cinf{M}$
satisfies that $\chi
\cc{z}^k|_U \in \AHol{U}$. We now use the above equation for
$h' = \chi \cc{z}^k$ and for $l=1$ yielding $0 =
\frac{2}{\im}g^{n\cc{k}}i_s(Z_n)\pi_z
\rWi^{(2)}|_U$ since $\tauWi(\chi \cc{z}^k)^{(1)}|_U = \d \cc{z}^k|_U$.
Consequently the only non-vanishing terms in $\pi_z
\rWi^{(2)}|_U$ can only be of symmetric degree $0$ as $g$ is
non-degenerate. But from the recursion formula for $\rWi$ and the
fact that $\sWi$ contains no term of symmetric degree $1$ it is
obvious that $\rWi$ contains no part of symmetric degree $0$ and
hence we have $\pi_z \rWi^{(2)}|_U =0$. Repeating this argument
with all holomorphic charts of $M$ this finally implies $\pi_z
\rWi^{(2)}=0$. We now assume that $\pi_z \rWi^{(s)} =0$ for $2\leq s
\leq m$ then the equation derived above yields for $l=m$ that
$0= \frac{1}{\nu} \pi_z ((\pi_z \rWi^{(m+1)})\circWi \tauWi
(h')^{(1)})|_U$. Repeating the above argument this now implies
$\pi_z \rWi^{(m+1)}=0$ and hence $\pi_z \rWi=0$ by induction. For
the proof of ii.) one proceeds analogously using the functions
$\chi z^k$ instead of $\chi \cc{z}^k$.
\end{PROOF}

\begin{LEMMA}
Let $U \subseteq M$ be an open subset of $M$.
\begin{enumerate}
\item
If we have $h' \astWi g|_U = h' g|_U$ for all $g,h'\in
\Cinf{M}$ with $h'|_U\in\AHol{U}$ then $\pi_z
\tauWi(h')|_U = h'|_U$ for all $h'\in \Cinf{M}$ with
$h'|_U\in\AHol{U}$.
\item
If we have $f \astWi h|_U = f h|_U$ for all $f,h\in
\Cinf{M}$ with $h|_U\in\Hol{U}$ then $\pi_{\cc{z}}\tauWi(h)|_U
= h|_U$ for all $h\in \Cinf{M}$ with $h|_U\in\Hol{U}$.
\end{enumerate}
\end{LEMMA}
\begin{PROOF}
It is easy to compute that $h' \astWi g = h' g +
\sum_{s=1}^\infty \sum_{r=1}^{2s-1} \sigma((\pi_z\tauWi(h')^{(r)})
\circWi (\pi_{\cc{z}}\tauWi(g)^{(2s-r)}))$. Since the summand
corresponding to the summation index $s$ is of $\nu$-degree $s$ the
condition $h' \astWi g|_U = h' g|_U$ implies that we have
$\sum_{r=1}^{2s-1} \sigma((\pi_z\tauWi(h')^{(r)})
\circWi (\pi_{\cc{z}}\tauWi(g)^{(2s-r)}))|_U=0$ for all $s\geq 1$.
We first consider the case $s=1$. Now let $(z,V)$ be a holomorphic
chart of $M$ and let $\chi \in \Cinf{M}$ with
$\mathrm{supp}(\chi)\subseteq V$ and $\chi|_{U'}=1$ for some open
subset $U'\subseteq V$ with $U \cap U' \neq \emptyset$. Then we
consider the functions $\chi \cc{z}^k  \in \Cinf{M}$ and from
$\tauWi(\chi \cc{z}^k)^{(1)}|_{U\cap U'}= \d\cc{z}^k|_{U\cap U'}$
we get from the above equation using $g= \chi \cc{z}^k$ that
$0=\frac{2\nu}{\im}
g^{n\cc{k}}i_s(Z_n)\pi_z\tauWi(h')^{(1)}|_{U\cap U'}$. As the
symmetric degree of $\pi_z\tauWi(h')^{(1)}$ is $1$ and $g$ is
non-degenerate we conclude that $\pi_z\tauWi(h')^{(1)}|_{U\cap
U'}=0$. Repeating this argument with all holomorphic charts such
that $U\cap U'\neq \emptyset$ we finally obtain $\pi_z
\tauWi(h')^{(1)}|_U =0$. Now we want to use induction on the total
degree to show that the assumption actually implies
$\pi_z\tauWi(h')^{(r)}|_U =0$ for all $r\geq 1$. To this end we
need the following:
\begin{SUBLEMMA}
Let $(z,V)$ be a holomorphic chart of $M$ around $p\in M$ with
$z(p)=0$. Further let $\chi \in \Cinf{M}$ with
$\mathrm{supp}(\chi)\subseteq V$ and $\chi|_{U'} =1$, where $p \in
U'\subseteq V$. Moreover, let $\astWi$ satisfy $h' \astWi g|_U = h'
g|_U$ for all $U\subseteq M$ and all $g, h'\in \Cinf{M}$ with $h'
\in \AHol{U}$. Then we have for all $r\geq 1$ that
$\tauWi(\chi\cc{z}^{k_1}\ldots \cc{z}^{k_r})^{(s)}|_p=0$ for $0\leq
s < r$ and $\tauWi(\chi\cc{z}^{k_1}\ldots \cc{z}^{k_r})^{(r)}|_p =
\d\cc{z}^{k_1}\vee\ldots \vee \d\cc{z}^{k_r}|_p$.
\end{SUBLEMMA}
\begin{INNERPROOF}
As the functions $\chi \cc{z}^k$ are locally anti-holomorphic on
$U'$ and as $\chi^r|_{U'}=1$ for $r \in \mathbb N$ we have
$\chi\cc{z}^{k_1}\ldots \cc{z}^{k_r}|_{U'}
= (\chi\cc{z}^{k_1})\ldots (\chi\cc{z}^{k_r})|_{U'} =
(\chi\cc{z}^{k_1})\astWi \ldots \astWi (\chi\cc{z}^{k_r})|_{U'}$.
To this equation we apply $\tauWi$ and consider the term of total
degree $s$ and get
\[
\tauWi(\chi\cc{z}^{k_1}\ldots
\cc{z}^{k_r})^{(s)}|_{U'} = \sum_{l_1 + \ldots + l_r=s}
\tauWi(\chi\cc{z}^{k_1})^{(l_1)}\circWi \ldots \circWi
\tauWi(\chi\cc{z}^{k_r})^{(l_r)}|_{U'}.
\]
For $s < r$ at least one $l_i$ is equal to $0$ and we get the first
assertion since $z(p)=0$ and $\tauWi(\chi \cc{z}^{k_i})^{(0)}|_{U'}
= \cc{z}^{k_i}|_{U'}$. For $s=r$ the only summand that does not vanish at
$p$ is the one for $l_1
=\ldots = l_r=1$ and observing that $\tauWi(\chi\cc{z}^{k_i})|_{U'}
= \d \cc{z}^{k_i}|_{U'}$ the explicit shape of $\circWi$ implies
the second statement of the sublemma.
\end{INNERPROOF}

We now assume that $\pi_z \tauWi(h')^{(r)}|_U=0$ for $1 \leq r \leq
m-1$, where $m\geq 2$. Using the equation $\sum_{r=1}^{2s-1}
\sigma((\pi_z\tauWi(h')^{(r)}) \circWi (\pi_{\cc{z}}
\tauWi(g)^{(2s-r)}))|_U=0$ for all $s$ with $2s \geq m+1$
this induction hypothesis implies that $0=\sum_{r=m}^{2s
-1}\sigma((\pi_z\tauWi(h')^{(r)})
\circWi (\pi_{\cc{z}}\tauWi(g)^{(2s-r)}))|_U$. Now let $(z,V)$
be a chart as in the sublemma then we use this equation for $g =
\chi
\cc{z}^{k_1} \ldots \cc{z}^{k_{2s-m}}$ at the point $p\in M$ and
obtain
\[
0 = \left(\frac{2\nu}{\im}\right)^{2s-m}g^{l_1\cc{k}_1}\ldots
g^{l_{2s-m}\cc{k}_{2s-m}} \sigma (i_s(Z_{l_1})\ldots
i_s(Z_{l_{2s-m}})\pi_z \tauWi(h')^{(m)})|_p
\]
for all $s$ with $2s \geq m+1$. Considering the equation for $m$
even and $m$ odd separately it is easy to see that this implies
$\pi_z \tauWi(h')^{(m)}|_p =0$. But as $p\in U$ was arbitrary we
may conclude that $\pi_z \tauWi(h')|_U=0$ proving i.). For the
proof of ii.) one proceeds quite analogously.
\end{PROOF}

Altogether we now have shown:

\begin{THEOREM}\label{WickCharThm}
Let $\DWi$ be the Fedosov derivation constructed for $\circWi$ from
the data $(\OWi,\sWi)$, where $\sWi$ contains no terms of symmetric
degree $1$ and let $\astWi$ be the corresponding star product on a
pseudo-K\"ahler manifold $(M,\omega,I)$ then we have the following
equivalences:
\begin{enumerate}
\item
\BEQAR{AeiquivallinksEq}
\begin{array}{l}
\pi_z\sWi=0 \\
\pi_z(1\otimes \OWi)=0
\end{array}
&\iff&
\pi_z \rWi = 0 \iff
\begin{array}{l}
\pi_z\tauWi (h')|_U= h'|_U \quad\forall h'
\in\Cinf{M}\\
\textrm{with } h'|_U \in \AHol{U}
\end{array}
\\
&\iff& h'\astWi g|_U = h' g|_U\quad \forall g,h'\in\Cinf{M}\textrm{
with } h'|_U \in \AHol{U},\label{AeiquivalfinlinksEq}
\EEQAR
\item
\BEQAR{AeiquivalrechtsEq}
\begin{array}{l}
\pi_{\cc{z}}\sWi=0 \\
\pi_{\cc{z}}(1\otimes \OWi)=0
\end{array}
&\iff& \pi_{\cc{z}} \rWi = 0 \iff
\begin{array}{l}
\pi_{\cc{z}}\tauWi (h)|_U=
h|_U \quad\forall h \in\Cinf{M}\\
\textrm{with } h|_U \in \Hol{U}
\end{array}
\\
&\iff& f \astWi h|_U = fh|_U\quad \forall f,h\in\Cinf{M}\textrm{
with } h|_U \in \Hol{U}.\label{AeiquivalfinrechtsEq}
\EEQAR
\item
Hence the star product $\astWi$ on $(M,\omega,I)$ is of Wick type
if and only if $\OWi$ is of type $(1,1)$ and $\pi_z \sWi =
\pi_{\cc{z}} \sWi =0$.
\end{enumerate}
\end{THEOREM}

After having uniquely characterized the star products $\astWi$ that
are of Wick type it is almost trivial to characterize the star
products $\astAW$ that are of anti-Wick type.

\begin{COROLLARY}
Let $\DAW$ be the Fedosov derivation obtained for $\circAW$ from
the data $(\OAW,\sAW)$, where $\sAW$ contains no terms of symmetric
degree $1$ and let $\astAW$ be the corresponding star product on a
pseudo-K\"ahler manifold $(M,\omega,I)$ then we have the following
equivalences:
\begin{enumerate}
\item
\BEQAR{AeiquivallinksAWEq}
\begin{array}{l}
\pi_z\sAW=0 \\
\pi_z(1\otimes \OAW)=0
\end{array}
&\iff&
\pi_z \rAW = 0 \iff
\begin{array}{l}
\pi_z\tauAW (h')|_U= h'|_U \quad\forall h'
\in\Cinf{M}\\
\textrm{with } h'|_U \in \AHol{U}
\end{array}
\\
&\iff& g\astAW h' |_U = g h' |_U\quad \forall
g,h'\in\Cinf{M}\textrm{ with } h'|_U \in
\AHol{U},\label{AeiquivalfinlinksAWEq}
\EEQAR
\item
\BEQAR{AeiquivalrechtsAWEq}
\begin{array}{l}
\pi_{\cc{z}}\sAW=0 \\
\pi_{\cc{z}}(1\otimes \OAW)=0
\end{array}
&\iff& \pi_{\cc{z}} \rAW = 0 \iff
\begin{array}{l}
\pi_{\cc{z}}\tauAW (h)|_U=
h|_U \quad\forall h \in\Cinf{M}\\
\textrm{with } h|_U \in \Hol{U}
\end{array}
\\
&\iff&  h\astAW f |_U = hf |_U\quad \forall f,h\in\Cinf{M}\textrm{
with } h|_U \in \Hol{U}.\label{AeiquivalfinrechtsAWEq}
\EEQAR
\item
Hence the star product $\astAW$ on $(M,\omega,I)$ is of anti-Wick
type if and only if $\OAW$ is of type $(1,1)$ and $\pi_z \sAW =
\pi_{\cc{z}} \sAW =0$.
\end{enumerate}
\end{COROLLARY}
\begin{PROOF}
We consider the star product $\astAW$ constructed from the data
$(\OAW,\sAW)$ then we have from Lemma \ref{AWWiRelLem} that the
star product $\astWi$ constructed from $(\OWi =\Pa \OAW,\sWi = \Pa
\sAW)$ satisfies $f \astAW g = \Pa ((\Pa g)\astWi (\Pa g))$ for all
$f,g\in \Cinf{M}[[\nu]]$. Moreover, in this case $\rWi = \Pa \rAW$
and $\tauWi (f) = \Pa \tauAW(P f)$ for all $f \in \Cinf{M}[[\nu]]$.
Using these relations the proof of the corollary is trivial using
the corresponding statements of Theorem \ref{WickCharThm}.
\end{PROOF}

After having uniquely characterized the star products $\astWi$
resp. $\astAW$ that are of Wick resp. anti-Wick type by some
conditions on the data $(\OWi,\sWi)$ resp. $(\OAW,\sAW)$ we can
show another interesting property of these star products namely
that they actually do not depend on $\sWi$ resp. $\sAW$. To this
end we consider two Fedosov star products $\astWi$ and
$\ast'_\TinyWi$ of Wick type that are obtained from the same formal
series of closed two-forms $\OWi = \Omega'_\TinyWi$ of type $(1,1)$
but different elements $\sWi$ and $s'_\TinyWi$. From Proposition
\ref{CharKlaProp} we already know that these star products are
equivalent and we even have an explicit construction for an
equivalence transformation. Letting $h_\TinyWi\in \mathcal W$ be
the element determined by the equation (\ref{hdiamRecEq}), where
$\diamond = \TinyWi$ and $C_\TinyWi =0$ we have that
$A_{h_\TinyWi}$ defined by $A_{h_\TinyWi} f:= \sigma\left(
\exp\left(\frac{1}{\nu}\adWi(h_\TinyWi)\right)\tauWi(f)\right)$ is
an equivalence from $\astWi$ to $\ast'_\TinyWi$ and our aim is to
show that the fact that both star products are of Wick type implies
that this equivalence transformation is equal to the identity.

\begin{LEMMA}
Let $\astWi$, $\ast'_\TinyWi$ and $h_\TinyWi$ be given as above
then we have
\BEQ{hProjEq}
\pi_z h_\TinyWi = \pi_{\cc{z}} h_\TinyWi =0.
\EEQ
\end{LEMMA}
\begin{PROOF}
We want to prove the statement by induction on the total degree and
hence we explicitly write the term of total degree $k\geq 3$ of
equation (\ref{hdiamRecEq}) as
\begin{eqnarray*}
{h_\TinyWi}^{(k)} &=&\delta^{-1}\left(\nabla {h_\TinyWi}^{(k-1)} -
\frac{1}{\nu}\sum_{\stackrel{k_1+l=k+1, l \geq 2}
{\scriptscriptstyle 3 \leq k_1 \leq k-1}}
\adWi(\rWi^{(l)}){h_\TinyWi}^{(k_1)}\right.\\
& &- \left.\sum_{r=0}^\infty
\frac{1}{r!}B_r \left(\frac{1}{\nu}\right)^r
\sum_{\stackrel{k_1+\ldots +k_r +l= k-1+2r}{\scriptscriptstyle
3\leq k_j\leq k- 1,l\geq 2}}
\adWi({h_\TinyWi}^{(k_1)})\ldots \adWi({h_\TinyWi}^{(k_r)})
(r'_\TinyWi -\rWi)^{(l)} \right),
\end{eqnarray*}
where $B_r$ denotes the $r^{\rm th}$ Bernoulli number that arises
from the Taylor expansion of $\frac{\frac{1}{\nu}\adWi
(h_\TinyWi)}{\exp\left(\frac{1}{\nu}\adWi (h_\TinyWi)\right)-\id}$
which is given by $\sum_{r=0}^\infty \frac{1}{r!} B_r
\left(\frac{1}{\nu}\right)^r\adWi(h_\TinyWi)^r$. For $k=3$
this yields using $B_0=1$ that ${h_\TinyWi}^{(3)} =
\delta^{-1}((\rWi - {r'}_\TinyWi)^{(2)})$ and hence we get
using Lemma \ref{ZerlegLem} and $\pi_z\rWi = \pi_z
r'_\TinyWi=\pi_{\cc{z}}\rWi=\pi_{\cc{z}}r'_\TinyWi=0$ which holds
since both star products are of Wick type that $\pi_z
{h_\TinyWi}^{(3)} = \pi_{\cc{z}} {h_\TinyWi}^{(3)}=0$. Let us now
assume that $\pi_z {h_\TinyWi}^{(m)}=0$ for $3\leq m \leq k-1$,
where $k \geq 4$ and consider
\begin{eqnarray*}
\pi_z {h_\TinyWi}^{(k)}\!\!&=&\!\!{\delta_z}^{\!\!\!-1}
\left(\nabla_z \pi_z {h_\TinyWi}^{(k-1)} -
\frac{1}{\nu}\sum_{\stackrel{k_1+l=k+1, l \geq 2}
{\scriptscriptstyle 3 \leq k_1 \leq k-1}}
\pi_z(\adWi(\rWi^{(l)}){h_\TinyWi}^{(k_1)})
\right.\\
\!\! & &\!\!-\left. \sum_{r=0}^\infty
\frac{1}{r!} B_r \left(\frac{1}{\nu}\right)^r
\!\!\!\sum_{\stackrel{k_1+\ldots +k_r +l= k-1+2r}{\scriptscriptstyle
3\leq k_j\leq k- 1, l\geq 2}}\!\!\!\!
\pi_z(\adWi({h_\TinyWi}^{(k_1)})\ldots \adWi({h_\TinyWi}^{(k_r)})
(r'_\TinyWi -\rWi)^{(l)}) \right).
\end{eqnarray*}
From the induction hypothesis we have ${\delta_z}^{\!\!\!-1}
(\nabla_z \pi_z {h_\TinyWi}^{(k-1)})=0$. Moreover, we have using
equation (\ref{picircVertrEq}) that
$\pi_z(\adWi(\rWi^{(l)}){h_\TinyWi}^{(k_1)})=
\pi_z((\pi_z\rWi^{(l)})
\circWi {h_\TinyWi}^{(k_1)}- (\pi_z {h_\TinyWi}^{(k_1)})\circWi
\rWi^{(l)})$. Here the first summand vanishes since $\astWi$ is of
Wick type whereas the second one vanishes by induction observing
that $3\leq k_1\leq k-1$. Now the contribution
$\pi_z(\adWi({h_\TinyWi}^{(k_1)})\ldots
\adWi({h_\TinyWi}^{(k_r)}) (r'_\TinyWi -\rWi)^{(l)})$ consists of
three types of terms. In case $r=0$ this contributes
${\delta_z}^{\!\!\!-1}\pi_z ((\rWi - r'_\TinyWi)^{(k-1)})$ which is
$0$ since both star products are of Wick type. For $r\neq 0$ the
whole expression consists of term of the shape $\pi_z
({h_\TinyWi}^{(k_j)}\circWi A )= \pi_z ((\pi_z
{h_\TinyWi}^{(k_j)})\circWi A )$ and
$\pi_z((r'_\TinyWi-\rWi)^{(l)}\circWi B) =
\pi_z((\pi_z (r'_\TinyWi-\rWi)^{(l)})\circWi B)$,
with certain elements $A\in \WL^1$ and $B\in \mathcal W$. Again
observing the range of summation the terms of the first kind vanish
by the induction hypothesis and the others vanish since
$\pi_z\rWi=\pi_z r'_\TinyWi=0$ proving $\pi_z {h_\TinyWi}^{(k)}=0$
and hence $\pi_z h_\TinyWi=0$ by induction. For the proof of
$\pi_{\cc{z}}h_\TinyWi$ one proceeds analogously.
\end{PROOF}

Using this result we can show that the star products $\astWi$ and
$\ast'_\TinyWi$ are not only equivalent but even equal.

\begin{PROPOSITION}\label{WiWistrichCoinProp}
Let $\astWi$ and $\ast'_\TinyWi$ be two star products of Wick type
on $(M,\omega,I)$ constructed from the data $(\OWi,\sWi)$ and
$(\OWi,s'_\TinyWi)$. Then the equivalence transformation
$A_{h_\TinyWi}$ from $\astWi$ to $\ast'_\TinyWi$ is equal to the
identity and hence the star products $\astWi$ and $\ast'_\TinyWi$
coincide.
\end{PROPOSITION}
\begin{PROOF}
We have $A_{h_\TinyWi} f= f + \sigma (\sum_{r=1}^\infty\frac{1}{r!}
\left(\frac{1}{\nu}\right)^r\adWi(h_\TinyWi)^r\tauWi(f))$ for all
$f\in\Cinf{M}[[\nu]]$ . Now the terms occurring in the sum are of
the shape
\[
\sigma(h_\TinyWi \circWi A) =\sigma((\pi_z h_\TinyWi)
\circWi A) \quad \textrm{ and } \quad \sigma(B \circWi h_\TinyWi)
= \sigma(B \circWi (\pi_{\cc{z}}h_\TinyWi))
\]
with certain elements $A,B \in \mathcal W$. From the preceding
lemma we have that both types of terms vanish and hence
$A_{h_\TinyWi}=\id$.
\end{PROOF}

\begin{COROLLARY}
The two star products of anti-Wick type $\astAW$ and
$\ast'_\TinyAW$ on $(M,\omega,I)$ constructed from the data
$(\OAW,\sAW)$ and $(\OAW,s'_\TinyAW)$ coincide.
\end{COROLLARY}
\begin{PROOF}
Considering the star products $\astWi$ and $\ast'_\TinyWi$ of Wick
type obtained from $(\OWi= \Pa \OAW, \sWi = \Pa \sAW)$ and
$(\OWi=\Pa \OAW, s'_\TinyWi = \Pa s'_\TinyAW)$ we have from
Proposition \ref{WiWistrichCoinProp} that $\astWi = \ast'_\TinyWi$.
But then the relation (\ref{astAWastWiPaOppEq}) implies that
$\astAW = \ast'_\TinyAW$.
\end{PROOF}

One should observe that the star products $\astWi$ resp. $\astAW$
that are not of Wick resp. anti-Wick type actually do depend on
$\sWi$ resp. $\sAW$. But in view of the above proposition and the
corollary we can always choose the simplest normalization condition
$\sWi=0$ resp. $\sAW=0$ when considering star products $\astWi$
resp. $\astAW$ of Wick resp. anti-Wick type. Hence the Fedosov
construction induces mappings
\BEQAR{OWiMapWick}
\nu Z^2_{\mbox{\rm\tiny dR}}(M,\mathbb C)^{1,1}[[\nu]]
\ni \OWi & \mapsto & \astWi \in
\{
\textrm{ star products of Wick type on $(M,\omega,I)$ }
\}\\
\label{OAWMapAntiWick}
\nu Z^2_{\mbox{\rm\tiny dR}}(M,\mathbb C)^{1,1}[[\nu]]
\ni \OAW & \mapsto & \astAW \in
\{
\textrm{ star products of anti-Wick type on $(M,\omega,I)$ }
\},
\EEQAR
where $\nu Z^2_{\mbox{\rm\tiny dR}}(M,\mathbb
C)^{1,1}[[\nu]]=\{\Omega \in
\nu \Ginf{\mbox{$\bigwedge$}^2T^*M}[[\nu]]\,|\, \d \Omega=0,
\pi^{1,1}\Omega =\Omega\}$ that shall turn out to be bijections. In
order to prove this we shall translate some of the results of
Karabegov obtained in \cite{Kar96} to the situation of star
products of Wick type in the next section.
\section{Identification of Star Products of Wick Type using
Karabegov's Method}
\label{KaraCharSec}
In this section we shall briefly report some general facts on star
products of Wick type that shall enable us to identify the Fedosov
star products of Wick type. Moreover, we give a direct elementary
proof (independent of a construction of such star products) of the
fact that star products of Wick type can be uniquely characterized
by some local equations.

Let $\starWi$ be a star product of Wick type on $(M,\omega,I)$ and
let $(z,V)$ denote a local holomorphic chart of $M$ such that
$V\subseteq M$ is a contractible open subset of $M$. It has been
shown in \cite[Prop. 1]{Kar96} that there are locally defined
formal functions $u_k \in \Cinf{V}[[\nu]]$ for $1 \leq k\leq
\dim_{\mathbb C}(M)$ such that
\BEQ{KarukEq}
u_k \starWi z^l - z^l \starWi u_k = - \nu \delta^l_k.
\EEQ
Observe that we are using a different sign convention for the
Poisson bracket causing the opposite sign in the above formula
compared to the one in \cite{Kar96}. Analogously one has that there
are locally defined formal functions $\cc{v}_l\in \Cinf{V}[[\nu]]$
for $1 \leq l \leq \dim_{\mathbb C}(M)$ such that
\BEQ{KarvcclEq}
\cc{v}_l \starWi \cc{z}^k - \cc{z}^k \starWi \cc{v}_l = \nu
\delta^{\cc{k}}_{\cc{l}}.
\EEQ
Using the fact that $\starWi$ is a star product of Wick type one
can show (cf. \cite[Lemma 2]{Kar96}) that the local functions
$u_k,\cc{v}_l\in \Cinf{V}[[\nu]]$ satisfy the equations
\BEQ{LinkRechWickMulKarEq}
f \starWi u_k = f u_k + \nu Z_k(f)\qquad\textrm{ and } \qquad
\cc{v}_l \starWi f = \cc{v}_l f + \nu \cc{Z}_l(f)
\EEQ
for all $f\in \Cinf{V}[[\nu]]$. Moreover, Karabegov considers
locally defined formal series of one-forms $\alpha,\beta\in
\Ginf{T^*V}[[\nu]]$ given by $\alpha:= - u_k \d z^k$ which is of type
$(1,0)$ and $\beta:= \cc{v}_l \d
\cc{z}^l$ which is of type $(0,1)$. As the
$\starWi$-right-multiplication with $u_k$ obviously commutes with
the $\starWi$-left-multiplication with $\cc{v}_l$ one in addition
obtains from (\ref{LinkRechWickMulKarEq}) that
$\cc{\partial}\alpha= \partial \beta$. Moreover, one can show that
this procedure yields a formal series of closed two-forms of type
$(1,1)$ on $M$ that does not depend on any of the choices made. In
the following the so-defined formal two-form that can be associated
to any star product $\starWi$ of Wick type shall be denoted by
$K(\starWi)$ and is referred to as Karabegov's characterizing form.
It is easy to see from the very definition that $K(\starWi) \in
\omega + \nu Z^2_{\mbox{\rm\tiny dR}}(M,\mathbb C)^{1,1}[[\nu]]$
and hence from the $\partial$-$\cc{\partial}$-Poincar\'{e} lemmas
one has that there exist formal local functions
$\varphi\in\Cinf{V}[[\nu]]$ on every contractible domain $V$ of
holomorphic coordinates such that $K(\starWi)|_V =
\partial\cc{\partial} \varphi$ and $\varphi$ is called a formal
local K\"ahler potential of $K(\starWi)$. Writing $\varphi =
\varphi_0 + \varphi_+$, where $\varphi_+ \in \nu \Cinf{V}[[\nu]]$
one evidently has that $\varphi_0$ is a local K\"ahler potential
for the pseudo-K\"ahler form $\omega$. With such a formal local
K\"ahler potential the equations
\BEQ{KahlPotStMultEq}
f \starWi Z_k(\varphi) = f Z_k(\varphi) + \nu Z_k(f)\qquad\textrm{
and }
\qquad
\cc{Z}_l(\varphi) \starWi f = \cc{Z}_l(\varphi) f + \nu \cc{Z}_l(f)
\EEQ
hold for all $f\in \Cinf{V}[[\nu]]$. Altogether one has the
following theorem:

\begin{THEOREM}{\rm\bf (\cite[Thm. 1]{Kar96})}\label{KStarSatz}
Let $\starWi$ be a star product of Wick type on a pseudo-K\"ahler
manifold $(M,\omega,I)$. Then $K(\starWi) \in
\omega + \nu Z^2_{\mbox{\rm\tiny dR}}(M,\mathbb C)^{1,1}[[\nu]]$
associates a formal series of closed two-forms of type $(1,1)$ on
$M$ which is a deformation of the pseudo-K\"ahler form $\omega$ to
this star product. In case $\varphi \in \Cinf{V}[[\nu]]$ is a
formal local K\"ahler potential of $K(\starWi)$ we have the
equations (\ref{KahlPotStMultEq}) for all $f\in \Cinf{V}[[\nu]]$.
\end{THEOREM}

Conversely in \cite[Sect. 4]{Kar96} Karabegov has shown that to
each such form $K$ as in the preceding theorem one can assign a
star product of Wick type such that the characterizing form of this
star product actually coincides with this given $K$. To this end
Karabegov has given an explicit construction of such a star product
extensively using local considerations. We now want to show that a
star product of Wick type in fact is completely determined by its
characterizing form which shall turn out to be the key result for
the proof that one can construct all star products of Wick type on
a pseudo-K\"ahler manifold using Fedosov's method with a suitably
chosen $\OWi$.

\begin{THEOREM}\label{WickTypLokCharSatz}
Let $\starWi$ and $\star'_\TinyWi$ be two star products of Wick
type on a pseudo-K\"ahler manifold $(M,\omega,I)$ and let
$\{\mathcal U_\gamma\}_{\gamma \in J}$ be a good open cover of $M$,
where the $\mathcal U_\gamma$ are the domains of holomorphic
coordinates of $M$. Moreover, let $K \in \omega +\nu
Z^2_{\mbox{\rm\tiny dR}}(M,\mathbb C)^{1,1}[[\nu]]$ and let
$\varphi_\gamma \in
\Cinf{\mathcal U_\gamma}[[\nu]]$ be formal local K\"ahler potentials
of $K$, i.e. $K|_{\mathcal U_\gamma} =
\partial\cc{\partial}\varphi_\gamma$.
\begin{enumerate}
\item
If we have for all $\gamma \in J$ and respectively all $f\in
\Cinf{\mathcal U_\gamma}[[\nu]]$ that
\BEQ{RechtCharEq}
f \starWi Z_k(\varphi_\gamma) = f \star'_\TinyWi
Z_k(\varphi_\gamma)
\EEQ
then $\starWi$ and $\star'_\TinyWi$ coincide.
\item
If we have for all $\gamma \in J$ and respectively all $f\in
\Cinf{\mathcal U_\gamma}[[\nu]]$ that
\BEQ{LinkCharEq}
\cc{Z}_l (\varphi_\gamma)\starWi f = \cc{Z}_l(\varphi_\gamma)
\star'_\TinyWi f
\EEQ
then $\starWi$ and $\star'_\TinyWi$ coincide.
\end{enumerate}
\end{THEOREM}
\begin{PROOF}
We assume that the bidifferential operators $C_i$ and $C'_i$
describing the star products $\starWi$ and $\star'_\TinyWi$
coincide for $0\leq i \leq k-1$, where $k\geq 1$ and want to show
that this implies that $C_k = C'_k$ proving by induction that the
star products coincide since all star products coincide in zeroth
order of $\nu$. From the associativity of both star products at the
order $k$ of $\nu$ one gets using the induction hypothesis that
$C_k - C'_k$ is a Hochschild cocycle. According to the
Hochschild-Kostant-Rosenberg-Theorem we thus have that there is a
differential operator $B$ on $\Cinf{M}$ and a two-form $A$ on $M$
such that $(C_k - C'_k)(f,g) = (\deltaH B)(f,g)+ A(X_f,X_g)$ for
all $f,g\in \Cinf{M}$. Here $X_f$ denotes the Hamiltonian vector
field corresponding to the function $f$, $\deltaH$ denotes the
Hochschild differential and as $C_k$ and $C'_k$ vanish on constants
$\deltaH B$ also vanishes on constants. Taking the anti-symmetric
part of the obtained equation and using that $\deltaH B$ is
symmetric we get
\[
A(X_f,X_g) = \frac{1}{2}\left(C_k(f,g)- C_k(g,f) - C_k'(f,g)+
C_k'(g,f)\right).
\]
Now let $z$ denote a holomorphic chart on $\mathcal U_\gamma$.
Since the bidifferential operator $(f,g)\mapsto A(X_f,X_g)$ is
obviously one-differential it is enough to evaluate it on the
coordinate functions $z^m$ and $\cc{z}^l$ to determine $A$. Using
the fact that $\starWi$ and $\star'_\TinyWi$ are of Wick type we
evidently have $C_k(z^l,z^m)= C'_k(z^l,
z^m)=C_k(\cc{z}^l,\cc{z}^m)=C'_k(\cc{z}^l,\cc{z}^m)=0$ and using
$X_{z^m}= -\frac{2}{\im}g^{m\cc{n}}\cc{Z}_n$ and $X_{\cc{z}^l}=
\frac{2}{\im}g^{n\cc{l}}Z_n$ we get
\[
0=A(X_{z^m},X_{z^l})= - 4 g^{m\cc{n}}g^{l\cc{j}}
A(\cc{Z}_n,\cc{Z}_j)\quad \textrm{ and }\quad
0=A(X_{\cc{z}^m},X_{\cc{z}^l}) = -4
g^{n\cc{m}}g^{j\cc{l}}A(Z_n,Z_j).
\]
Repeating the same argument for all $\gamma \in J$ and a
corresponding local holomorphic chart this yields with the fact
that $g$ is non-degenerate that $A$ is of type $(1,1)$ and in a
local chart we write $A= A_{i\cc{j}} \d z^i \wedge \d\cc{z}^j$ such
that $A(X_f,X_g)= 4 A_{i\cc{j}}g^{i \cc{l}}g^{n\cc{j}} (
\cc{Z}_l(f)Z_n(g)-\cc{Z}_l(g) Z_n(f))$. To proceed we also need some
special property of the symmetric part $\deltaH B$ of $C_k - C'_k$.
\begin{SUBLEMMA}
Let $\deltaH B$ be a bidifferential operator vanishing on constants
that in local holomorphic coordinates $(z,\mathcal U_\gamma)$ has
the shape
\BEQ{deltaHBlocEq1}
(\deltaH B)(f,g) =
\sum_{K,\cc{L}} C^{K;\cc{L}}\left(\frac{\partial^{|K|}f}{\partial
z^K}\frac{\partial^{|\cc{L}|}g}{\partial
\cc{z}^L}+\frac{\partial^{|K|}g}{\partial
z^K}\frac{\partial^{|\cc{L}|}f}{\partial \cc{z}^L}
\right)
\EEQ
with local functions $C^{K;\cc{L}}\in \Cinf{\mathcal U_\gamma}$.
Then in local holomorphic coordinates $\deltaH B$ has the shape
\BEQ{deltaHBlocEq2}
(\deltaH B)(f,g) = C^{k\cc{l}}(Z_k(f)\cc{Z}_l(g)+
Z_k(g)\cc{Z}_l(f)),
\EEQ
where $C^{k\cc{l}}$ denotes the components of a tensor field $C
\in \Ginf{TM\otimes TM}$ of type $(1,1)$.
\end{SUBLEMMA}
\begin{INNERPROOF}
From $\deltaH^2=0$ we get $0=f (\deltaH B)(g,h)- (\deltaH B)(fg,h)+
(\deltaH B)(f,gh)-(\deltaH B)(f,g)h$ for all $f,g,h\in \Cinf{M}$.
Evaluating this equation for $g|_U, h|_U \in \Hol{U}$, where $U
\subseteq M$ is an open subset with $U\cap \mathcal U_\gamma \neq
\emptyset$ one obtains with the above shape of $\deltaH B$ that
\[
0 = \sum_{K,\cc{L}}C^{K;\cc{L}}\left(
\frac{\partial^{|K|}(gh)}{\partial
z^K}\frac{\partial^{|\cc{L}|}f}{\partial \cc{z}^L}
- g \frac{\partial^{|K|}h}{\partial
z^K}\frac{\partial^{|\cc{L}|}f}{\partial
\cc{z}^L}- h\frac{\partial^{|K|}g}{\partial
z^K}\frac{\partial^{|\cc{L}|}f}{\partial
\cc{z}^L}
\right)
\]
on $U\cap \mathcal U_\gamma$. Here one has to observe that $|K|$
and $|\cc{L}|$ are at least one since $\deltaH B$ vanishes on
constants. From this we may conclude that
$\sum_{K,\cc{L}}C^{K;\cc{L}}\frac{\partial^{|K|}f}{\partial
z^K}\frac{\partial^{|\cc{L}|}g}{\partial
\cc{z}^L}=\sum_{k,\cc{L}}C^{k;\cc{L}}Z_k(f)
\frac{\partial^{|\cc{L}|}g}{\partial
\cc{z}^L}$ on $\mathcal U_\gamma$. Analogously one finds
$\sum_{K,\cc{L}}C^{K;\cc{L}}
\frac{\partial^{|K|}f}{\partial
z^K}\frac{\partial^{|\cc{L}|}g}{\partial
\cc{z}^L}= \sum_{K,\cc{l}} C^{K;\cc{l}}
\frac{\partial^{|K|}f}{\partial
z^K}\cc{Z}_l(g)$ repeating the above argument with $f,g,h\in
\Cinf{M}$ such that $f|_U,g|_U \in \AHol{U}$. Altogether the
proven local equations prove the statement since it is easy to see
that the local functions $C^{i\cc{j}}$ define a tensor field.
\end{INNERPROOF}

From the fact that both star products are of Wick type it is
obvious that $(\deltaH B)(f,g)=\frac{1}{2}((C_k - C_k')(f,g)+(C_k -
C_k')(g,f))$ is of the shape as in the sublemma yielding
\[
(C_k - C_k')(f,g)=C^{k\cc{l}}(Z_k(f)\cc{Z}_l(g)+ Z_k(g)\cc{Z}_l(f))
+4
A_{i\cc{j}} g^{i\cc{l}}g^{n\cc{j}}(\cc{Z}_l(f)Z_n(g)-
\cc{Z}_l(g)Z_n(f)).
\]
Again using that the star products are of Wick type this implies
$C^{n\cc{l}} = -4A_{i\cc{j}} g^{i\cc{l}}g^{n\cc{j}}$ and finally
\[
(C_k - C_k')(f,g)= -8 A_{i\cc{j}}
g^{i\cc{l}}g^{n\cc{j}}Z_n(f)\cc{Z}_l(g).
\]
But now we need the conditions (\ref{RechtCharEq}) resp.
(\ref{LinkCharEq}) to proceed with the proof. Considering
(\ref{RechtCharEq}) in order $k$ of the formal parameter we get
from the induction hypothesis that $C_k(f, Z_n(\varphi_{\gamma
0}))-C'_k(f,Z_n(\varphi_{\gamma 0}))=0$ for $f\in \Cinf{\mathcal
U_\gamma}$, where $\varphi_{\gamma 0}$ denotes a K\"ahler potential
for $\omega$ on $\mathcal U_\gamma$. Using the expression for $C_k
-C'_k$ derived above this yields $0 = -
4\im A_{n\cc{j}}g^{m\cc{j}}Z_m(f)$ and hence $A_{n\cc{j}}=0$. As
$\gamma \in J$ was arbitrary we get $A=0$ proving i.). The final
step to prove ii.) is analogous.
\end{PROOF}

The statement of this theorem means that a star product of Wick
type is uniquely determined by one of the equations
(\ref{KahlPotStMultEq}) in case it is valid in every holomorphic
chart $(z,V)$ and hence it is determined by $K$.

\begin{REMARK}\label{AntiWickRem}
All the statements about star products of Wick type made in this
section can be easily transferred to star products of anti-Wick
type. For completeness we give the important relations but omit the
proofs since they are easily done using the relation between star
products of Wick type and those of anti-Wick type explained in
Section \ref{NotDefSec}. For a star product $\starAW$ one defines
the characterizing form $\cc{K}(\starAW)\in
\omega + \nu Z^2_{\mbox{\rm\tiny dR}}(M,\mathbb C)^{1,1}[[\nu]]$ by
$\cc{K}(\starAW):= \Pa K({(\starAW)}_{\pa,\opp})$. Denoting by
$\cc{\varphi}\in
\Cinf{V}[[\nu]]$ a formal K\"ahler potential of $\cc{K}(\starAW)$
on $V$ that is obtained from a formal K\"ahler potential $\varphi$
of $K({(\starAW)}_{\pa,\opp})$ by $\cc{\varphi}=\Pa \varphi$ one
has
\BEQ{KahlPotStMultAntiEq}
Z_k(\cc{\varphi})\starAW f  = Z_k(\cc{\varphi}) f  - \nu
Z_k(f)\qquad\textrm{ and }
\qquad
f \starAW \cc{Z}_l(\cc{\varphi})  = f \cc{Z}_l(\cc{\varphi}) - \nu
\cc{Z}_l(f),
\EEQ
where $(z,V)$ is a holomorphic chart and $f \in \Cinf{V}[[\nu]]$.
As in Theorem \ref{WickTypLokCharSatz} the form $\cc{K}(\starAW)$
uniquely determines the star product of anti-Wick type by one of
these local equations valid in every holomorphic chart $(z,V)$ of
$M$.
\end{REMARK}
\section{Universality of Fedosov's
Construction for Star Products of Wick Type}
\label{FedUnivSec}
In this section we want to compute Karabegov's characterizing form
$K(\astWi)$ of the Fedosov star products of Wick type constructed
in this paper. Furthermore, our result shall enable us to prove
that the mapping defined in equation (\ref{OWiMapWick}) is a
bijection. To determine $K(\astWi)$ we need according to Section
\ref{KaraCharSec} local formal functions $u_k \in \Cinf{V}[[\nu]]$
with
\[
u_k \astWi z^l - z^l \astWi u_k = - \nu \delta^l_k,
\]
where $(z,V)$ is a local holomorphic chart of $M$. Denoting by
$\varphi_0$ a local K\"ahler potential of $\omega$ on $V$ we
obviously have
\[
- \delta^{l}_k = \{Z_k(\varphi_0), z^l\} = - \sigma(\Lie_{Z_k}
\tauWi(z^l))
\]
since the Lie derivative commutes with $\sigma$. Defining
$u_{k0}:=Z_k(\varphi_0)$ we thus have to show the existence of
local formal functions $u_{k+} \in \nu \Cinf{V}[[\nu]]$ such that
\[
-\sigma(\Lie_{Z_k} \tauWi(z^l)) = \frac{1}{\nu}
\sigma(\adWi(\tauWi(u_{k0}+u_{k+}))\tauWi(z^l)).
\]
For the further computations we need the following formula:

\begin{LEMMA}\label{LieDerLem}
For all vector fields $X\in \Ginf{TM}$ the Lie derivative $\Lie_X
:\WL \to \WL$ can be expressed in the following way:
\BEQAR{LieallgEq}
\Lie_X &=& \DWi i_a(X) + i_a(X)\DWi +
\frac{1}{\nu}\adWi(i_a(X)\rWi)
+ i_s(X)\nonumber\\ & & + (\d z^n\otimes 1) i_s(\nabla_{Z_n}X)+
(\d\cc{z}^m
\otimes 1) i_s(\nabla_{\cc{Z}_m}X).
\EEQAR
\end{LEMMA}
\begin{PROOF}
The proof of this formula is the same as for the analogous formula
using $\circF$ and $\DF$ instead of $\circWi$ and $\DWi$ and can be
found in \cite[Appx. A]{Neu99}.
\end{PROOF}

Using the preceding lemma for $X=Z_k$ we get
\[
- \sigma(\Lie_{Z_k}\tauWi(z^l)) = -\sigma\left(\frac{1}{\nu}
\adWi(i_a(Z_k)\rWi)\tauWi(z^l)+ i_s(Z_k)\tauWi(z^l)\right),
\]
where we have used that $\nabla_{\cc{Z}_m} Z_k=0$ and the obvious
equations $\DWi \tauWi(z^l) =0 =i_a(Z_k) \tauWi(z^l)$ and $\sigma(
(\d z^n\otimes 1) i_s(\nabla_{Z_n}Z_k)\tauWi(z^l)) =0$. Defining
$a_k
:= i_{Z_k}\omega \otimes 1 \in \mathcal W(V)$ one easily finds
$i_s(Z_k) = -\frac{1}{\nu} \adWi(a_k)$ and we get
\begin{eqnarray*}
- \sigma(\Lie_{Z_k}\tauWi(z^l))
&=&-\frac{1}{\nu}
\sigma((\pi_z(i_a(Z_k)\rWi - a_k))\circWi (\pi_{\cc{z}}\tauWi(z^l)))
\\
& &+ \frac{1}{\nu}\sigma((\pi_z\tauWi(z^l))\circWi
(\pi_{\cc{z}}(i_a(Z_k)\rWi - a_k))),
\end{eqnarray*}
where we have used $\sigma= \pi_z\pi_{\cc{z}}$ and
(\ref{picircVertrEq}). Obviously we have $\pi_z a_k =0$ and
$\pi_{\cc{z}}a_k = a_k$ since $\omega$ is of type $(1,1)$ and
$\pi_z i_a(Z_k)\rWi=i_a(Z_k)\pi_z \rWi =0$ from Theorem
\ref{WickCharThm} such that we finally obtain
\BEQ{deltaKomplEq}
-\nu \delta_k^l =\sigma((\pi_z\tauWi(z^l))\circWi
(\pi_{\cc{z}}i_a(Z_k)\rWi - a_k)).
\EEQ
From this equation it is obvious that $\pi_{\cc{z}}i_a(Z_k)\rWi -
a_k$ has to be strongly related to the Fedosov-Taylor series of the
function $u_k$ that is to be found and we have:

\begin{LEMMA}\label{tauukBedLem}
Let $\varphi_0$ be a local K\"ahler potential of $\omega$ on the
open contractible domain $V$ of a holomorphic chart $z$ of $M$ and
let $u_{k0} \in \Cinf{V}$ be defined by $u_{k0}:= Z_k(\varphi_0)$.
In case there exist local formal functions $u_{k+}\in \nu
\Cinf{V}[[\nu]]$ such that we have
\BEQ{tauukTayEq}
\pi_{\cc{z}} \tauWi (u_k) = u_k +  i_{Z_k}
\omega \otimes 1-\pi_{\cc{z}}i_a(Z_k)\rWi
\EEQ
for $u_k = u_{k0}+ u_{k+}$ then $u_k \astWi z^l- z^l\astWi u_k
=-\nu \delta^l_k$.
\end{LEMMA}
\begin{PROOF}
For the proof we compute using equation (\ref{deltaKomplEq})
\[
-\nu\delta_k^l = \sigma((\pi_z\tauWi(z^l)) \circWi(u_k -
\pi_{\cc{z}}\tauWi(u_k)) )= u_k z^l - z^l\astWi u_k= u_k \astWi z^l
- z^l \astWi u_k,
\]
where the last equality follows from the fact that $\astWi$ is of
Wick type.
\end{PROOF}

So it remains to show that there exist such local functions as in
(\ref{tauukTayEq}) and to get more information about their concrete
shape.

\begin{PROPOSITION}\label{tauukExGestProp}
With the notation of Lemma \ref{tauukBedLem} one has the following
statements:
\begin{enumerate}
\item
For all $u_{k+}\in \nu\Cinf{V}[[\nu]]$ one has $\sigma(u_k +
i_{Z_k} \omega \otimes 1-\pi_{\cc{z}}i_a(Z_k)\rWi ) = u_k$ and
\BEQ{DWiccZtauccZEq}
-\delta_{\cc{z}}A_k +\nabla_{\cc{z}} A_k -
\frac{1}{\nu}\pi_{\cc{z}}(\rWi \circWi
A_k ) = 1\otimes (\cc{\partial}u_{k+}- i_{Z_k}\OWi),
\EEQ
where $A_k:=u_k + i_{Z_k} \omega
\otimes 1-\pi_{\cc{z}}i_a(Z_k)\rWi \in \mathcal W(V)$.
\item
One has $\pi_{\cc{z}} \tauWi (u_k) = u_k + i_{Z_k}
\omega \otimes 1 - \pi_{\cc{z}}i_a(Z_k)\rWi$
if and only if $\cc{\partial}u_{k+}- i_{Z_k}\OWi=0$. In case
$u_{k+}= Z_k(\varphi_+)$, where $\varphi_+\in\nu\Cinf{V}[[\nu]]$ is
a local formal K\"ahler potential for $\OWi$, i.e.
$\partial\cc{\partial}\varphi_+ =
\OWi|_V$ one has $\pi_{\cc{z}}\tauWi(u_k)=u_k +
i_{Z_k}\omega\otimes 1 - \pi_{\cc{z}}i_a(Z_k)\rWi$.
\item
Karabegov's characterizing form $K(\astWi)$ of the star product
$\astWi$ of Wick type is given by
\BEQ{KastWiEq1}
K(\astWi)= \omega + \OWi.
\EEQ
\item
For all $f\in\Cinf{V}[[\nu]]$ and $u_k$ as in ii.) one has
\BEQ{rechtsukMultEq}
f \astWi u_k = f u_k + \nu Z_k(f).
\EEQ
\end{enumerate}
\end{PROPOSITION}
\begin{PROOF}
Obviously we have $\sigma(u_k + i_{Z_k} \omega \otimes
1-\pi_{\cc{z}}i_a(Z_k)\rWi)=u_k$ for all $u_{k+}\in
\Cinf{V}[[\nu]]$ since $\sigma =
\pi_{\cc{z}}\pi_z$ and $\pi_z(i_a(Z_k)\rWi)=0$. Further we compute
for the proof of i.)
\[
\delta_{\cc{z}}(u_k + i_{Z_k}\omega\otimes 1 -
\pi_{\cc{z}}i_a(Z_k)\rWi)= 1\otimes i_{Z_k}\omega - \pi_{\cc{z}}
\delta i_a(Z_k)\rWi=1\otimes i_{Z_k}\omega -
\pi_{\cc{z}}(i_s(Z_k)\rWi - i_a(Z_k)\delta \rWi),
\]
where we have used $\delta i_a(Z_k)+i_a(Z_k)\delta = i_s(Z_k)$ and
(\ref{KomRelholantiholEq}). With the equation (\ref{littlerEq}) for
$\rWi$ this yields using $\pi_{\cc{z}}i_a(Z_k)R=0$,
$\pi_{\cc{z}}\rWi=0$ and
\begin{eqnarray*}
\lefteqn{\pi_{\cc{z}}(i_a(Z_k)(\rWi \circWi \rWi))}\\
&=&
\pi_{\cc{z}}((i_a(Z_k)\rWi)\circWi (\pi_{\cc{z}}\rWi)) -
\pi_{\cc{z}}(\rWi \circWi (i_a(Z_k)\rWi))=
-\pi_{\cc{z}}(\rWi \circWi (i_a(Z_k)\rWi))
\end{eqnarray*}
the intermediate result
\begin{eqnarray*}
\lefteqn{\delta_{\cc{z}}(u_k + i_{Z_k}\omega\otimes 1 -
\pi_{\cc{z}}i_a(Z_k)\rWi)}\\
&=&1\otimes i_{Z_k}\omega -
\pi_{\cc{z}}i_s(Z_k)\rWi + \pi_{\cc{z}}(i_a(Z_k)\nabla\rWi)+
\frac{1}{\nu} \pi_{\cc{z}}(\rWi \circWi (i_a(Z_k)\rWi)) +
1\otimes i_{Z_k} \OWi,
\end{eqnarray*}
since $\OWi$ is of type $(1,1)$. Further we have from
$\nabla_{\cc{Z}_l}Z_k=0$ and $\nabla_{\cc{Z}_l}\omega=0$ that
$\nabla_{\cc{Z}_l}i_{Z_k}\omega$ vanishes also, such that we obtain
using $\nabla i_a(Z_k)+i_a(Z_k)\nabla = \Lie_{Z_k}- (\d z^n\otimes
1)i_s(\nabla_{Z_n}Z_k)$ and (\ref{KomRelholantiholEq})
\begin{eqnarray*}
\lefteqn{\nabla_{\cc{z}}(u_k + i_{Z_k}\omega \otimes 1 -
\pi_{\cc{z}} i_a(Z_k)\rWi)}\\ &=& 1\otimes \cc{\partial} u_k -
\pi_{\cc{z}}(\Lie_{Z_k}\rWi -
(\d z^n\otimes 1)i_s(\nabla_{Z_n}Z_k)\rWi -
i_a(Z_k)\nabla\rWi)=1\otimes \cc{\partial} u_k+
\pi_{\cc{z}}i_a(Z_k)\nabla\rWi
\end{eqnarray*}
since $\pi_{\cc{z}}\rWi=0$ implies $\pi_{\cc{z}}\Lie_{Z_k}\rWi=0$
also. Finally we compute again with $\pi_{\cc{z}}\rWi=0$
\[
- \frac{1}{\nu}\pi_{\cc{z}}(\rWi \circWi (u_k +
i_{Z_k}\omega\otimes 1- \pi_{\cc{z}}i_a(Z_k)\rWi)) = -
\pi_{\cc{z}} i_s(Z_k)\rWi + \frac{1}{\nu}\pi_{\cc{z}}(\rWi\circWi
(i_a(Z_k)\rWi)),
\]
where besides (\ref{picircVertrEq}) we have used the trivial but
helpful identity $0=\pi_{\cc{z}}((i_{Z_k}\omega\otimes 1)\circWi
\rWi)$ and $i_s(Z_k)= -\frac{1}{\nu}\adWi
( i_{Z_k}\omega\otimes 1 )$. Combining our three intermediate
results one finds (\ref{DWiccZtauccZEq}) since $u_{k0}=
Z_k(\varphi_0)$ and hence $\cc{\partial}u_{k0}=i_{Z_k}\omega$
proving i.). According to Proposition \ref{tauWickeinfachProp} ii.)
we thus have that $A_k=\pi_{\cc{z}}\tauWi(u_k)$ if and only if
$\cc{\partial}u_{k+}- i_{Z_k}\OWi=0$. But this equation is
obviously satisfied in case one has $u_{k+}= Z_k (\varphi_+)$ with
a local K\"ahler potential $\varphi_+$ of $\OWi$. According to
Lemma \ref{tauukBedLem} we then have $u_k
\astWi z^l -z^l \astWi u_k = -\nu \delta_k^l$ and we get using the
definition of the characterizing form $K(\astWi)$ with
$u_k=Z_k(\varphi)=Z_k(\varphi_0+\varphi_+)$ that
\[
K(\astWi)= - \cc{\partial}(u_k \d z^k)=- \cc{\partial}(u_k)\wedge
\d z^k= \d z^k \wedge i_{Z_k}(\omega+\OWi)=\omega + \OWi
\]
proving iii.). For the proof of iv.) we compute for
$f\in\Cinf{V}[[\nu]]$ using $\pi_{\cc{z}}\tauWi(u_k) = u_k
+i_{Z_k}\omega \otimes 1 - \pi_{\cc{z}}i_a(Z_k)\rWi$,
$\DWi\tauWi(f)=0$ and (\ref{picircVertrEq})
\begin{eqnarray*}
f \astWi u_k &=& \sigma( \tauWi(f) \circWi (u_k +i_{Z_k}\omega
\otimes 1 - i_a(Z_k)\rWi))\\
&=&f u_k +
\frac{2\nu}{\im}g^{n\cc{m}}\sigma(i_s(Z_n)\tauWi(f))(i_{Z_k}\omega)
(\cc{Z}_m)+ \sigma(i_a(Z_k) \adWi (\rWi)\tauWi(f))\\ & &-
\sigma((i_a(Z_k)\rWi)\circWi \tauWi(f))\\
&=&f u_k + \nu \sigma(i_s(Z_k)\tauWi(f)) + \nu
\sigma(i_a(Z_k)(\nabla -\delta)\tauWi(f))-
\sigma((\pi_z i_a(Z_k)\rWi)\circWi \tauWi(f)).
\end{eqnarray*}
Further we obtain by reason of $\pi_z i_a(Z_k)\rWi
=i_a(Z_k)\pi_z \rWi=0$ and $i_a(Z_k)(\nabla-\delta)\tauWi(f)=
(\nabla_{Z_k}-i_s(Z_k))\tauWi(f)$ that
\[
f \astWi u_k =fu_k + \nu \sigma(\nabla_{Z_k}\tauWi(f))= fu_k + \nu
Z_k (f),
\]
and the proposition is proven.
\end{PROOF}

Slightly modifying the proofs of Lemma \ref{tauukBedLem} and
Proposition \ref{tauukExGestProp} one can also show the existence
of locally defined formal functions $\cc{v}_l \in \Cinf{V}[[\nu]]$
such that $\cc{v}_l \astWi \cc{z}^k - \cc{z}^k\astWi \cc{v}_l=\nu
\delta^{\cc{k}}_{\cc{l}}$. For completeness we just give the results
and omit the proof since by Theorem \ref{WickTypLokCharSatz} it is
enough to know the properties of one of the functions $u_k$ resp.
$\cc{v}_l$ to be able to identify the star product $\astWi$.
Moreover, the proofs given above can easily be transferred to the
following statements.

\begin{LEMMA}\label{tauccvlBedLem}
Let $\varphi_0$ be a local K\"ahler potential of $\omega$ on the
open contractible domain $V$ of a holomorphic chart $z$ of $M$ and
let $\cc{v}_{l0} \in \Cinf{V}$ be defined by $\cc{v}_{l0}:=
\cc{Z}_l(\varphi_0)$. In case there exist local formal functions
$\cc{v}_{l+}\in
\nu \Cinf{V}[[\nu]]$ such that we have
\BEQ{tauccvlTayEq}
\pi_z \tauWi (\cc{v}_l) = \cc{v}_l -  i_{\cc{Z}_l}
\omega \otimes 1 + \pi_z i_a(\cc{Z}_l)\rWi
\EEQ
for $\cc{v}_l = \cc{v}_{l0}+ \cc{v}_{l+}$ then $\cc{v}_l \astWi
\cc{z}^k- \cc{z}^k \astWi \cc{v}_l= \nu \delta_{\cc{l}}^{\cc{k}}$.
\end{LEMMA}
\begin{PROPOSITION}\label{tauccvlExGestProp}
With the notation of Lemma \ref{tauccvlBedLem} one has the
following statements:
\begin{enumerate}
\item
For all $\cc{v}_{l+}\in \nu\Cinf{V}[[\nu]]$ one has
$\sigma(\cc{v}_l - i_{\cc{Z}_l}\omega \otimes 1+\pi_z
i_a(\cc{Z}_l)\rWi ) = \cc{v}_l$ and
\BEQ{DWiZtauZEq}
-\delta_z\cc{B}_l +\nabla_z \cc{B}_l +
\frac{1}{\nu}\pi_z(\cc{B}_l \circWi \rWi ) =
1\otimes (\partial \cc{v}_{l+} + i_{\cc{Z}_l}\OWi),
\EEQ
where $\cc{B}_l:=\cc{v}_l - i_{\cc{Z}_l} \omega
\otimes 1 + \pi_z i_a(\cc{Z}_l)\rWi \in \mathcal W(V)$.
\item
One has $\pi_z \tauWi (\cc{v}_l) = \cc{v}_l - i_{\cc{Z}_l}
\omega \otimes 1 + \pi_z i_a(\cc{Z}_l)\rWi$ if and only if
$\partial\cc{v}_{l+}+  i_{\cc{Z}_l}\OWi=0$. In case $\cc{v}_{l+}=
\cc{Z}_l(\varphi_+)$, where $\varphi_+\in \nu \Cinf{V}[[\nu]]$ is a
local formal K\"ahler potential for $\OWi$, i.e.
$\partial\cc{\partial}\varphi_+ =
\OWi|_V$ one has $\pi_z\tauWi(\cc{v}_l)=\cc{v}_l -
i_{\cc{Z}_l}\omega\otimes 1 + \pi_z i_a(\cc{Z}_l)\rWi$.
\item
Here also Karabegov's characterizing form $K(\astWi)$ of the star
product $\astWi$ of Wick type is given by
\BEQ{KastWiEq2}
K(\astWi)= \omega + \OWi.
\EEQ
\item
For all $f\in\Cinf{V}[[\nu]]$ and $\cc{v}_l$ as in ii.) one has
\BEQ{linksccvlMultEq}
\cc{v}_l \astWi f = \cc{v}_lf  + \nu \cc{Z}_l(f).
\EEQ
\end{enumerate}
\end{PROPOSITION}

The assertions iii.) and iv.) of Propositions \ref{tauukExGestProp}
and \ref{tauccvlExGestProp} represent the generalizations of the
theorem proven by Karabegov in \cite[Sect. 4]{Kar99} for the
special case $\OWi=0$. In the following corollary we just collect
the more or less trivial analogous statements for the star products
$\astAW$ of anti-Wick type.

\begin{COROLLARY}\label{KantiWickCor}
Let $\astAW$ be the Fedosov star product of anti-Wick type
constructed from $\OAW\in \nu Z^2_{\mbox{\rm\tiny dR}}(M,\mathbb
C)^{1,1}[[\nu]]$ then its characterizing form is given by
\BEQ{antiWickCharFormEq}
\cc{K}(\astAW) = \omega + \OAW.
\EEQ
Denoting by $\cc{\varphi}$ a local formal K\"ahler potential of
$\cc{K}(\astAW)$, i.e. $\partial\cc{\partial}\cc{\varphi}
|_V=\cc{K}(\astAW)$ and defining $\cc{u}_l:=
\cc{Z}_l(\cc{\varphi}), v_k := Z_k(\cc{\varphi})\in
\Cinf{V}[[\nu]]$ we have
\BEQ{cculvkEq}
v_k\astAW f  = v_k  f  - \nu Z_k(f)\qquad\textrm{ and }
\qquad
f \astAW \cc{u}_l  = f \cc{u}_l - \nu
\cc{Z}_l(f),
\EEQ
where $(z,V)$ is a holomorphic chart and $f \in \Cinf{V}[[\nu]]$.
\end{COROLLARY}
\begin{PROOF}
We only have to prove the formula for $\cc{K}(\astAW)$ since the
equations (\ref{cculvkEq}) are known to be true from Remark
\ref{AntiWickRem} equation (\ref{KahlPotStMultAntiEq}). Using Lemma
\ref{AWWiRelLem} it is obvious that $(\astAW)_{\pa,\opp}$ is the
star product of Wick type obtained from $\OWi= \Pa \OAW$. Using
(\ref{KastWiEq1}) this implies the assertion since $\Pa^2=\id$.
\end{PROOF}

We now present the main result of our paper:

\begin{THEOREM}\label{FedWickUnivSatz}
{\rm\bf (Universality of Fedosov's construction for star products
of Wick type)} For all star products $\starWi$ of Wick type on a
pseudo-K\"ahler manifold $(M,\omega,I)$ there is a Fedosov
construction with $\circWi$ such that the obtained star product
$\astWi$ coincides with $\starWi$. Moreover, the mapping
\BEQ{OWimapastWi}
\nu Z^2_{\mbox{\rm\tiny dR}}(M,\mathbb C)^{1,1}[[\nu]]
\ni \OWi  \mapsto  \astWi \in
\{
\textrm{ star products of Wick type on $(M,\omega,I)$ }
\}
\EEQ
induced by the Fedosov construction is a bijection and its inverse
is given by $\astWi \mapsto K(\astWi)-\omega$, where $K(\astWi)$
denotes Karabegov's characterizing form of $\astWi$.
\end{THEOREM}
\begin{PROOF}
The universality of the Fedosov construction just means that the
mapping in equation (\ref{OWimapastWi}) is surjective and hence we
have to show that for any given star product $\starWi$ of Wick type
there is a $\OWi$ such that $\astWi=\starWi$. Given $\starWi$ we
consider $K(\starWi)$ according to Theorem \ref{KStarSatz} and have
that the equations (\ref{KahlPotStMultEq}) are fulfilled for a
local formal K\"ahler potential $\varphi$ of $K(\starWi)$ on
$V\subseteq M$ and $f \in \Cinf{V}[[\nu]]$. Using the Fedosov
construction for $\OWi = K(\starWi)-\omega$ we have from
Propositions \ref{tauukExGestProp} and \ref{tauccvlExGestProp} that
$f\astWi Z_k(\varphi) = f Z_k(\varphi) + \nu Z_k(f)$ and
$\cc{Z}_l(\varphi) \astWi f = \cc{Z}_l(\varphi) f + \nu
\cc{Z}_l(f)$. With Theorem \ref{WickTypLokCharSatz} each of these
equations implies $\astWi = \starWi$ proving surjectivity. From
$K(\astWi)=\omega + \OWi$ the fact that the mapping in
(\ref{OWimapastWi}) is injective and the shape of its inverse are
obvious.
\end{PROOF}

The preceding theorem can immediately be transferred to the star
products of anti-Wick type that can all be obtained from a Fedosov
construction using $\circAW$ and the mapping according to equation
(\ref{OAWMapAntiWick}) also turns out to be a bijection.

A quite remarkable consequence of this theorem is:

\begin{DEDUCTION}
All star products $\starWi$ ($\starAW$) of (anti-)Wick type on a
pseudo-K\"ahler manifold $(M,\omega,I)$ are of Vey type.
\end{DEDUCTION}
\begin{PROOF}
By the above theorem all star products of (anti-)Wick type can be
obtained by some Fedosov construction using $\circWi$ ($\circAW$).
But from Proposition \ref{VeyTypProp} we have that all star
products $\astWi$ ($\astAW$) even those that are not of (anti-)Wick
type are of Vey type.
\end{PROOF}

\begin{DEDUCTION}{\rm\bf (\cite[Thm. 3]{Kar98})} The characteristic
classes $c(\starWi)$ resp. $c(\starAW)$ of star products of Wick
resp. anti-Wick type on a pseudo-K\"ahler manifold $(M,\omega,I)$
are related to the characterizing forms $K(\starWi)$ resp.
$\cc{K}(\starAW)$ by $c(\starWi)=
\frac{[K(\starWi)]}{\nu}-\frac{[\rho]}{\im}$ resp. $c(\starAW)=
\frac{[\cc{K}(\starAW)]}{\nu}+ \frac{[\rho]}{\im}$, where $\rho$
denotes the Ricci form.
\end{DEDUCTION}
\begin{PROOF}
With the results of Proposition \ref{tauukExGestProp} iii.) resp.
Corollary \ref{KantiWickCor} and Proposition \ref{CharKlaProp} i.)
the statements are obviously true for the star products $\astWi$
resp. $\astAW$ of Wick resp. anti-Wick type. But as the Fedosov
constructions are universal for these types of star products they
in fact hold in general.
\end{PROOF}

Another application of our results is the unique characterization
of the Hermitian star products of Wick and of anti-Wick type that
is those star products that have the complex conjugation denoted by
$\C$ incorporated as an anti-automorphism (cf. \cite[Sect.
5]{Neu99}).

\begin{PROPOSITION}
A star product $\starWi$ ($\starAW$) of (anti-)Wick type on a
pseudo-K\"ahler manifold $(M,\omega,I)$ is Hermitian if and only if
the characterizing form $K(\starWi)$ ($\cc{K}(\starAW)$) is real,
i.e. $\C K(\starWi) = K(\starWi)$
($\C\cc{K}(\starAW)=\cc{K}(\starAW)$).
\end{PROPOSITION}
\begin{PROOF}
Let $\starWi$ be a Hermitian star product of Wick type and let $u_k
\in\Cinf{V}[[\nu]]$ such that $f \starWi u_k = f u_k + \nu Z_k(f)$ for
all $f\in\Cinf{V}[[\nu]]$. Using $\C \nu = -\nu$ since the formal
parameter is assumed to be purely imaginary this implies that
$\C(f)\C(u_k) - \nu \cc{Z}_k(\C f) = \C(f \astWi u_k) = (\C u_k)
\astWi (\C f)$ and hence we have $\cc{v}'_k \starWi g = \cc{v}'_k g
+ \nu \cc{Z}_k(g)$ for all $g\in \Cinf{V}[[\nu]]$, where
$\cc{v}'_k:=- \C u_k$. Using the definition of the characterizing
form and $\C^2=\id$ this yields
\[
K(\starWi)= \partial (\cc{v}'_k \d\cc{z}^k) =\C \C\partial C(-u_k
\d z^k) = \C \cc{\partial}(- u_k \d z^k) = \C K(\starWi).
\]
Now let $\starWi$ be a star product of Wick type such that $\C
K(\starWi)=K(\starWi)$ and consider a star product $\astWi=\starWi$
which exists according to Theorem \ref{FedWickUnivSatz}. Using
$K(\astWi)=\omega + \OWi$ we obviously have $\C\OWi =\OWi$. We now
want to show that this implies that $\astWi$ is Hermitian. To this
end we first observe that the fibrewise product $\circWi$ satisfies
$\C (a \circWi b)= (-1)^{kl}(\C b) \circWi (\C a)$ for $a\in
\WL^k,b\in \WL^l$. Using this it is straightforward to see that
$\C\rWi$ satisfies the same equations as $\rWi$ does yielding
$\C\rWi =\rWi$. But then we have $\DWi \C = \C \DWi$ and hence $\C
\tauWi(f)=\tauWi(\C f)$ for all $f\in \Cinf{M}[[\nu]]$ implying by
an easy computation that $\astWi$ and therefore $\starWi$ is
Hermitian. The proof of the analogous statements for the star
products of anti-Wick type is completely analogous.
\end{PROOF}

To conclude this section we want to establish a relation to the
star products with separation of variables considered by Karabegov.
To this end we have to drop our convention to use the purely
imaginary formal parameter $\nu$ and have to replace it by
$\im\lambda$ and consider the star product defined by
\BEQ{KaraStarDef}
f \astK g := g f +\sum_{l=1}^\infty (\im\lambda)^l C_l(g,f)
\EEQ
for $f,g \in \Cinf{M}$, where the bidifferential operators $C_l$
describe the star product $\astWi$ of Wick type by $f \astWi g = fg
+ \sum_{l=1}^\infty \nu^l C_l(f,g)$. Obviously $\astK$ is a star
product with separation of variables on $(M,-\omega, I)$ and using
Karabegov's original definition of the characterizing form (cf.
\cite[Sect. 3]{Kar96}) it is an easy task to compute $K(\astK) =
\omega + \OWi(\im \lambda)$ using Proposition
\ref{tauukExGestProp}. Moreover, it is easy to conclude from
Proposition \ref{CharKlaProp} and the properties of the
characteristic class (cf. \cite[Lemma 5.2 i.)]{Neu99} and
\cite[Thm. 6.4]{GR99}) that $c(\astK)=\frac{\im}{\lambda}
[K(\astK)]-\im [\rho]$. It stands to reason that one can also give
a Fedosov construction that directly yields the star product
$\astK$. We just give a sketch of the necessary modifications:
Consider $(\WL_\lambda,\circK)$, where $\WL_\lambda$ is the same
object as $\WL$, where $\nu$ has been replaced by $\lambda$ and
$\circK$ is defined by
\BEQ{circKDefEq}
a\circK b:= \mu \circ \exp \left( 2\lambda
g^{k\cc{l}}i_s(\cc{Z}_l)\otimes i_s(Z_k)\right)(a\otimes b).
\EEQ
For this fibrewise product one obtains a Fedosov derivation $\DK
:= -\delta + \nabla -\frac{1}{\im\lambda}\adK(\rK)$ with
$\DK^2=0$, where $\rK$ is the solution of the equations
\BEQ{rKEq}
\delta\rK = \nabla\rK - \frac{1}{\im\lambda} \rK \circK \rK - R -
1 \otimes \OWi(\im \lambda)\qquad\textrm{ and }\qquad
\delta^{-1}\rK = 0.
\EEQ
It is easy to see that $\rK = - \rWi(\im \lambda)$ and that
$\tauK(f) = \tauWi(f)(\im\lambda)$ for all $f\in \Cinf{M}$, but
this actually implies that the obtained star product coincides with
the one defined in equation (\ref{KaraStarDef}).
\section*{Outlook and Further Questions}
Let us conclude with a few remarks on this approach to star
products of Wick type. First we would like to point out that the
Fedosov construction has the advantage that one is only dealing
with global geometric objects and is not forced to use expressions
of bidifferential operators in local charts. Moreover, the Fedosov
setting allows for very detailed investigations of algebraic
properties of the obtained star products just using their
construction and there is no need to have explicit closed formulas
for the star products. So it might be useful to consider the
following questions and topics in this framework:
\begin{enumerate}
\item
It should be possible to study representations of the star products
of Wick type and of anti-Wick type within Fedosov's framework.
Moreover, the question of Morita equivalence of these star products
can be discussed in this setting by explicitly constructing
equivalence bimodules which is subject of the work \cite{NeuWal02}.
\item
With a more detailed understanding of the representation theory of
the considered star products it should be possible to establish a
relation to geometric quantization and to Berezin-Toeplitz
quantization starting from the formal framework (cf.
\cite{KarSch00}).
\item
For a star product $\starWi$ of Wick type one can define the formal
Berezin transform $B$ as the unique formal series of differential
operators on $\Cinf{M}$ such that for any open subset $U\subseteq
M$ and $f,g\in \Cinf{M}$ with $f|_U \in\Hol{U}, g|_U\in \AHol{U}$
the relation $B(f g)= f\starWi g$ holds. Then the star product
$\starAW$ defined by $f \starAW g :=B^{-1}((Bf)\starWi(Bg))$ turns
out to be a star product of anti-Wick type and it would be
interesting to find a relation between $K(\starWi)$ and
$\cc{K}(\starAW)$ in this case and to obtain a more explicit
description of $B$.
\item
Finally it should be possible to consider different ways of
obtaining star products (possibly of Wick resp. anti-Wick type) by
phase space reduction (for instance BRST quantization as in
\cite{BorHerWal00}, or the analogue of the procedure in
\cite{Fed98}) in case there is a group acting on $(M,\omega,I)$ not
only preserving the symplectic form but also the complex structure.
Again this will be discussed in a forthcoming project.
\end{enumerate}
\begin{appendix}
\section{Consequences of the Complex Structure}
\label{ComStrConSec}
Using the complex structure on $(M,\omega,I)$ one can define
several splittings of the mappings which are involved into the
Fedosov construction which were of great advantage in Sections
\ref{WickCharSec} and \ref{FedUnivSec}. First it is obvious from
$I^2=-\id$ that the bundle isomorphism $I$ of the complexified
tangent bundle $TM$ has the fibrewise eigenvalues $\pm \im$ and
that $TM = TM^{1,0}\oplus TM^{0,1}$, where $TM^{1,0}$ denotes the
bundle of the $+\im$-eigenspaces and $TM^{0,1}$ denotes the bundle
of $-\im$-eigenspaces. Moreover, all sections $X\in \Ginf{TM}$ can
be uniquely decomposed into a section of type $(1,0)$ and a section
of type $(0,1)$. Considering the dual bundles $T^*M^{1,0}$ and
$T^*M^{0,1}$ to $TM^{1,0}$ and $TM^{0,1}$ it is obvious that the
cotangent bundle can also be written as $T^*M = T^*M^{1,0}\oplus
T^*M^{0,1}$. This decomposition naturally gives rise to analogous
decompositions of sections $T\in \Ginf{\bigvee^s T^*M}$ and
$\beta\in
\Ginf{\bigwedge^{a'}T^*M}$ such that we have
\BEQ{wedgeZerEq}
\Ginf{\mbox{$\bigwedge$}T^*M} = \bigoplus_{a'=0}^{\dim_\mathbb R(M)}
\bigoplus_{p'+q'=a'} \Ginf{\mbox{$\bigwedge$}^{p',q'}T^*M}.
\EEQ
Here $\mbox{$\bigwedge$}^{p',q'}T^*M$ denotes the bundle with the
characteristic fibre $\mbox{$\bigwedge$}^{p',q'}T_x^*M$, $x\in M$
and $\mbox{$\bigwedge$}^{p',q'}T_x^*M$ denotes the subspace of
$\mbox{$\bigwedge$}^{a'}T_x^*M$ with $a'=p'+q'$ that is generated
by elements of the shape $v'\wedge w'$ with $v'\in
\mbox{$\bigwedge$}^{p'} T_x^*M^{1,0}$ and $w'
\in \mbox{$\bigwedge$}^{q'}T_x^*M^{0,1}$.
In addition this decomposition induces natural projections
\BEQ{wedgeProjEq}
\pi^{p',q'} : \Ginf{\mbox{$\bigwedge$}^{a'}T^*M} \to
\Ginf{\mbox{$\bigwedge$}^{p',q'}T^*M},  \textrm{ where } a'=p'+q'
\EEQ
onto the part of type $(p',q')$. Analogously we have
\BEQ{veeZerEq}
\mathsf{X}_{s=0}^\infty \Ginf{
\mbox{$\bigvee$}^s T^*M} = \mathsf{X}_{s=0}^\infty
\bigoplus_{p+q=s} \Ginf{\mbox{$\bigvee$}^{p,q}T^*M}.
\EEQ
Here $\mbox{$\bigvee$}^{p,q}T^*M$ denotes the bundle with the
characteristic fibre $\mbox{$\bigvee$}^{p,q}T_x^*M$, $x\in M$ and
$\mbox{$\bigvee$}^{p,q}T_x^*M$ denotes the subspace of
$\mbox{$\bigvee$}^sT_x^*M$ with $s=p+q$ that is generated by
elements of the shape $v\vee w$ with $v\in
\mbox{$\bigvee$}^{p} T_x^*M^{1,0}$ and $w \in
\mbox{$\bigvee$}^{q}T_x^*M^{0,1}$.
Again there are the induced projections
\BEQ{veeProjEq}
\pi^{p,q} : \Ginf{\mbox{$\bigvee$}^s T^*M} \to
\Ginf{\mbox{$\bigvee$}^{p,q}T^*M}, \textrm{ where } s= p+q
\EEQ
onto the sections of type $(p,q)$. The defined projections extend
in a natural way to $\Ginf{\bigwedge T^*M}$ and
$\mathsf{X}_{s=0}^\infty \Ginf{\bigvee^s T^*M}$ setting
$\pi^{p',q'}(\beta):=0$ for $\beta\in \Ginf{\bigwedge^{a'}T^*M}$
with $a'\neq p'+q'$ and $\pi^{p,q}(T):=0$ for $T\in\Ginf{\bigvee^s
T^*M}$ with $s\neq p+q$. Using these projections we define the
mappings $\pi_s^{p,q},\pi_a^{p',q'}:\WL \to \WL$ on factorized
sections $T\otimes \beta\in \WL$ by
\BEQ{piapisDefEq}
\pi_s^{p,q}(T\otimes \beta) := (\pi^{p,q}T )\otimes \beta\qquad
\textrm{ and }\qquad \pi_a^{p',q'}(T\otimes \beta):=
T \otimes(\pi^{p',q'}\beta).
\EEQ
From these mappings we obtain further projections onto the purely
holomorphic resp. purely anti-holomorphic part of the symmetric and
the anti-symmetric part of $\WL$ by
\BEQ{pizcczsymDefEq}
\pi_{s,z} := \sum_{p=0}^\infty\pi_s^{p,0} \qquad \textrm{ resp. }
\qquad\pi_{s,\cc{z}} := \sum_{q=0}^\infty\pi_s^{0,q}
\EEQ
and
\BEQ{pizcczantisymDefEq}
\pi_{a,z} := \sum_{p'=0}^{\dim_{\mathbb R}(M)}\pi_a^{p',0} \qquad
\textrm{ resp. } \qquad\pi_{a,\cc{z}} := \sum_{q'=0}^{\dim_{\mathbb R}(M)}
\pi_a^{0,q'}.
\EEQ
With these projections we moreover define the projections onto the
totally holomorphic resp. totally anti-holomorphic part of $\WL$ by
\BEQ{purholantiholProjDef}
\pi_z :=\pi_{s,z}\pi_{a,z}= \pi_{a,z}\pi_{s,z}\qquad \textrm{ resp. }
\qquad \pi_{\cc{z}} := \pi_{s,\cc{z}}\pi_{a,\cc{z}}=
\pi_{a,\cc{z}} \pi_{s,\cc{z}}.
\EEQ

\begin{LEMMA}\label{piLem}
The mappings defined in equations (\ref{pizcczsymDefEq}),
(\ref{pizcczantisymDefEq}) and (\ref{purholantiholProjDef}) are
projections, i.e. we have $\pi_{s,z}\pi_{s,z}=\pi_{s,z}$,
$\pi_{a,z}\pi_{a,z}=\pi_{a,z}$,
$\pi_{s,\cc{z}}\pi_{s,\cc{z}}=\pi_{s,\cc{z}}$,
$\pi_{a,\cc{z}}\pi_{a,\cc{z}}=\pi_{a,\cc{z}}$ and hence we also
have $\pi_z\pi_z=\pi_z$ and
$\pi_{\cc{z}}\pi_{\cc{z}}=\pi_{\cc{z}}$. Furthermore, these
projections are homomorphisms of the undeformed product on $\WL$.
In addition $\pi_{a,z}$ and $\pi_{a,\cc{z}}$ are homomorphisms of
the fibrewise Wick product $\circWi$. For all $a,b\in \WL$ we
moreover have the equations
\BEQ{pisymcircVertrEq}
\pi_{s,z}(a\circWi b) = \pi_{s,z}((\pi_{s,z} a)\circWi b)\qquad
\textrm{ and }\qquad \pi_{s,\cc{z}}(a\circWi b) = \pi_{s,\cc{z}}(
a \circWi (\pi_{s,\cc{z}}b))
\EEQ
and
\BEQ{picircVertrEq}
\pi_z(a\circWi b) = \pi_z((\pi_z a)\circWi b)\qquad
\textrm{ and }\qquad \pi_{\cc{z}}(a\circWi b) = \pi_{\cc{z}}(
a \circWi (\pi_{\cc{z}}b)).
\EEQ
In addition we have $\sigma = \pi_z
\pi_{\cc{z}}=\pi_{\cc{z}}\pi_z$.
\end{LEMMA}
\begin{PROOF}
The proof of this lemma is straightforward using the very
definitions and the explicit shape of $\circWi$.
\end{PROOF}

Now we consider the mappings $\delta$, $\delta^*$ and $\nabla$ more
closely. Using a local holomorphic chart of $M$ the mappings
$\delta_z,\delta_{\cc{z}},{\delta_z}^{\!\!\!*},
{\delta_{\cc{z}}}^{\!\!\!*},\nabla_z,\nabla_{\cc{z}}:\WL\to
\WL$ are obviously well-defined by
\BEQ{splitEqDef}
\begin{array}{rcl}
\delta_z &:=& (1\otimes \d z^k)i_s(Z_k),\\
\delta_{\cc{z}}&:=&(1\otimes \d \cc{z}^l) i_s(\cc{Z}_l),
\end{array}
\begin{array}{rcl}
{\delta_z}^{\!\!\!*}&:=&(\d z^k \otimes 1)i_a(Z_k),\\
{\delta_{\cc{z}}}^{\!\!\!*}&:=&(\d \cc{z}^l\otimes 1)i_a(\cc{Z}_l),
\end{array}
\begin{array}{rcl}
\nabla_z&:=&(1\otimes \d z^k)\nabla_{Z_k},\\
\nabla_{\cc{z}}&:=&(1\otimes \d \cc{z}^l) \nabla_{\cc{Z}_l}
\end{array}
\EEQ
and we have $\delta = \delta_z +
\delta_{\cc{z}}$, $\delta^*={\delta_z}^{\!\!\!*} +
{\delta_{\cc{z}}}^{\!\!\!*}$ and $\nabla = \nabla_z +
\nabla_{\cc{z}}$. Completely analogously to the definition of
$\delta^{-1}$ we now define ${\delta_z}^{\!\!\!-1}$ for $a\in \WL$
with $(\d z^i\otimes 1)i_s(Z_i)a = k a$ and $(1\otimes \d z^i)
i_a(Z_i)a=l a$ by
\BEQ{delzminus1Def}
{\delta_z}^{\!\!\!-1}a := \left\{
    \begin {array} {cl}
    \frac{1}{k+l} {\delta_z}^{\!\!\!*}a &
    \mbox { in case } k+l \neq 0 \\
    0 & \mbox { in case } k+l=0
    \end {array}
    \right..
\EEQ
For $a\in \WL$ with $(\d\cc{z}^i\otimes 1)i_s(\cc{Z}_i)a = k a$ and
$(1\otimes \d\cc{z}^i) i_a(\cc{Z}_i)a=l a$ we analogously define
\BEQ{delzbarminus1Def}
{\delta_{\cc{z}}}^{\!\!\!-1}a := \left\{
    \begin {array} {cl}
    \frac{1}{k+l} {\delta_{\cc{z}}}^{\!\!\!*}a &
    \mbox { in case } k+l \neq 0
    \\ 0 & \mbox { in case } k+l=0
    \end {array}
    \right..
\EEQ

In the following two lemmas we just collect some properties of the
defined splittings of the mappings $\delta$, $\delta^{-1}$ and
$\nabla$.

\begin{LEMMA}\label{ZerlegLem}
For all $a \in \WL$ we have the decompositions
\BEQ{HodgeholaholZerEq}
{\delta_z}^{\!\!\!-1} \delta_z a + \delta_z {\delta_z}^{\!\!\!-1} a
+
\pi_{\cc{z}}a=a\qquad
\textrm{ and }\qquad
{\delta_{\cc{z}}}^{\!\!\!-1} \delta_{\cc{z}} a + \delta_{\cc{z}}
{\delta_{\cc{z}}}^{\!\!\!-1} a +
\pi_z a=a.
\EEQ
Furthermore, we have the following relations:
\BEQ{KomRelholantiholEq}
\begin{array}{rcl}
\pi_z \delta &=& \delta_z \pi_z,\\
\pi_{\cc{z}} \delta &=& \delta_{\cc{z}} \pi_{\cc{z}},
\end{array}
\quad
\begin{array}{rcl}
\pi_z \nabla &=& \nabla_z \pi_z, \\
\pi_{\cc{z}} \nabla &=& \nabla_{\cc{z}} \pi_{\cc{z}},
\end{array}
\quad
\begin{array}{rcl}
\pi_z \delta^{-1} &=& {\delta_z}^{\!\!\!-1} \pi_z, \\
\pi_{\cc{z}} \delta^{-1} &=& {\delta_{\cc{z}}}^{\!\!\!-1}
\pi_{\cc{z}}.
\end{array}
\EEQ
\end{LEMMA}
\begin{PROOF}
Again the proof is straightforward using the definitions and the
property of the pseudo-K\"ahler connection to be compatible with
the complex structure.
\end{PROOF}

\begin{LEMMA}\label{SquareLem}
The holomorphic and anti-holomorphic parts of the mappings $\delta$
and $\nabla$ satisfy the following identities:
\BEQ{zcczIdentities}
\begin{array}{l}
\delta_z^2=\delta_{\cc{z}}^2=[\delta_z,\delta_{\cc{z}}]=0,\\
\nabla_z^2= \nabla_{\cc{z}}^2=0,
\quad [\nabla_z,\nabla_{\cc{z}}]=-\frac{1}{\nu} \adWi(R),
\end{array}
\qquad
\begin{array}{r}
[\delta_z,\nabla_z]=[\delta_{\cc{z}},\nabla_{\cc{z}}]=0,
\\\
[\delta_z,\nabla_{\cc{z}}] = [\delta_{\cc{z}},\nabla_z]=0.
\end{array}
\EEQ
For all $b \in \pi_z(\WL^k),c \in \pi_{\cc{z}}(\WL^l), a\in
\WL^m$ we have:
\BEQAR{DerzIdentities}
\begin{array}{rcl}
\delta_z \pi_z (b \circWi a)&=& \pi_z(\delta_z b\circWi a)+ (-1)^k
\pi_z(b \circWi \delta a), \\
\nabla_z \pi_z (b \circWi a)&=& \pi_z(\nabla_z b\circWi a)
+(-1)^k\pi_z(b \circWi \nabla a),
\end{array}
\\
\begin{array}{rcl}
\label{DercczIdentities}
\delta_{\cc{z}} \pi_{\cc{z}} (a\circWi c)&=& (-1)^m \pi_{\cc{z}}
(a \circWi \delta_{\cc{z}}c) + \pi_{\cc{z}}(\delta a \circWi c),\\
\nabla_{\cc{z}} \pi_{\cc{z}} (a\circWi c)&=& (-1)^m\pi_{\cc{z}} (a
\circWi \nabla_{\cc{z}}c) + \pi_{\cc{z}}(\nabla a \circWi c).
\end{array}
\EEQAR
\end{LEMMA}
\begin{PROOF}
The statements in equation (\ref{zcczIdentities}) directly follow
from the properties of $\delta$ and $\nabla$ considering the
equations $\delta^2=[\delta,\nabla]=0$ and
$\nabla^2=-\frac{1}{\nu}\adWi(R)$ sorted with respect to the
holomorphic and anti-holomorphic degrees. Using Lemma
\ref{ZerlegLem} it is easy to verify (\ref{DerzIdentities}) and
(\ref{DercczIdentities}) since $\delta$ and $\nabla$ are
super-derivations with respect to $\circWi$.
\end{PROOF}

\section{An Application of Banach's Fixed Point Theorem}
\label{BanFixPktSec}
In this appendix we collect some consequences of Banach's fixed
point theorem that were very useful for divers proofs in Sections
\ref{FedProdEquCharClaSec} and \ref{WickCharSec}. Our presentation
which is slightly more general than it is actually needed mainly
follows \cite[Appx. 1.2]{Wal99}.

First let us recall Banach's fixed point theorem in its usual
formulation. We consider a metric space $(\mathcal M,d)$ and a
mapping $L:\mathcal M\to \mathcal M$ from this space to itself,
then $L$ is called contracting in case there is a $q\in \mathbb R$
with $0\leq q < 1$ such that we have $d(L(x),L(y))\leq q d(x,y)$
for all $x,y \in \mathcal M$. In other words this means that $L$ is
Lipschitz continuous with Lipschitz constant $q$. In case
$(\mathcal M,d)$ is in addition complete we have that such a
contracting mapping $L$ has a unique fixed point according to
Banach's fixed point theorem:

\begin{LEMMA}\label{BanLem}
Let $(\mathcal M,d)$ be a complete metric space and let $L :
\mathcal M\to \mathcal M$ be a contracting mapping, then $L$ has a
unique fixed point $x_0\in \mathcal M$ that can be obtained by
iteration $x_0 = \lim_{n \to \infty} L^n(x)$, where $x  \in
\mathcal M$ is arbitrary.
\end{LEMMA}

The idea for the following application of this statement on the one
hand is to find an appropriate metric such that the considered
space becomes a complete metric space and on the other hand to
guarantee that the relevant mappings are contracting with respect
to this metric. Now let $\mathsf R$ denote a ring. Then we consider
the Cartesian product $V:= \mathsf{X}_{k=0}^\infty V_k$ of $\mathsf
R$-modules $V_k$ for $k\in \mathbb N$ which again is a $\mathsf
R$-module. We define the order of an element $v= (v_k)_{k\in
\mathbb N}\in V$ by $o(v):=\min \{k \in \mathbb N\,|\,v_k \neq 0\}$ for
$v\neq 0$ and $o(0):= + \infty$. Using this definition we define
the valuation $\varphi: V \to \mathbb Q$ by $\varphi(v):=
2^{-o(v)}$ for $v\neq 0$ and $\varphi(0):=0$ and get the following:

\begin{LEMMA}\label{VkomLem}
With the notations from above $(V= \mathsf{X}_{k=0}^\infty
V_k,d_\varphi)$, where $d_\varphi$ is defined by
$d_\varphi(v,w):=\varphi(v-w)$ for $v,w\in V$, is a complete
ultrametric space. Moreover, a mapping $L:V \to V$ is contracting
with respect to $d_\varphi$ if and only if there is a $0 < k
\in\mathbb N \cup \{+ \infty\}$ such that
\begin{equation}\label{conEq}
o(L(v)-L(w))\geq k + o(v-w)\qquad\forall v, w \in V.
\end{equation}
In this case $2^{-k}$ is a Lipschitz constant for $L$, where we
have set $2^{-\infty}:=0$.
\end{LEMMA}
\begin{PROOF}
The proof is just a slight modification of the proof in case all
the $\mathsf R$-modules $V_k$ coincide with $V_0$. The proof for
this situation, where $V= V_0[[\nu]]$ is the space of formal power
series with values in $V_0$ and the topology induced by $d_\varphi$
is the $\nu$-adic topology, can be found in several textbooks e.g.
\cite[p. 388]{Kas95} or in the article \cite[Prop. 2]{BorWal98}.
\end{PROOF}

Since the Fedosov algebra $\WL$ is equal to the Cartesian product
$\WL = {\mathsf X}_{k=0}^\infty (\WL)^{(k)}$, where $(\WL)^{(k)}$
denotes the subspace of the elements that are homogeneous of degree
$k$ with respect to the total degree, i.e. $\Deg a = k a$ for all
$a \in (\WL)^{(k)}$, we can apply the above lemma and obtain:

\begin{COROLLARY}\label{FedConAbbCor}
\begin{enumerate}
\item
A mapping $L : \WL \to \WL$ has a unique fixed point in $\WL$, in
case there is a $0 < k \in \mathbb N$ such that for all $a,b \in
\WL$ the total degree of the term with lowest total degree in $L(a)
- L(b)$ is at least higher by $k$ than the total degree of the term
of lowest total degree in $a-b$.
\item
If $L:\WL \to \WL$ raises the total degree at least by $1$ then $L$
is contracting.
\item
Let $\pi:\WL \to \WL$ denote a projection that is continuous with
respect to the topology induced by the above ultrametric. Then the
statements analogous to i.) and ii.) hold for a mapping $L :
\pi(\WL) \to \pi(\WL)$ and the unique fixed point lies in
$\pi(\WL)$.
\end{enumerate}
\end{COROLLARY}
\begin{PROOF}
The statements i.) and ii.) are obvious from the very definitions.
For the proof of iii.) one just has to observe that the continuity
of $\pi$ and the fact that $\pi$ is a projection imply that
$\pi(\WL)= \ker(\id - \pi)\subseteq \WL$ is closed and hence a
complete ultrametric space with respect to the restriction of the
metric on $\WL$ to $\pi(\WL)$.
\end{PROOF}

\end{appendix}


\begin{thebibliography}{99}
\bibitem{BayFla78}
{\sc Bayen, F., Flato, M., Fr\o nsdal, C., Lichnerowicz, A.,
Sternheimer, D.:} {\it Deformation Theory and Quantization.} Ann.
Phys. {\bf 111}, Part I: 61--110, Part II: 111--151 (1978).
%
\bibitem{BerCahGut97}
{\sc Bertelson, M., Cahen, M., Gutt, S.:} {\it Equivalence of star
products.} Class. Quant. Grav. {\bf 14}, A93--A107 (1997).
%
\bibitem{BorHerWal00}
{\sc Bordemann, M., Herbig, H.-C., Waldmann, S.:} {\it BRST
Cohomology and Phase Space Reduction in Deformation Quantization.}
Commun. Math. Phys. {\bf 210}, 107--144 (2000).
%
\bibitem{BorNeu98a}
{\sc Bordemann, M., Neumaier, N., Waldmann, S.:} {\it Homogeneous
Fedosov Star Products on Cotangent Bundles I: Weyl and Standard
Ordering with Differential Operator Representation.} Commun. Math.
Phys. {\bf 198}, 363--396 (1998).
%
\bibitem{BorWal96a}
{\sc Bordemann, M., Waldmann, S.:} {\it A Fedosov Star Product of
Wick Type for K\"{a}hler Manifolds.} Lett. Math. Phys. {\bf 41},
243--253 (1997).
%
\bibitem{BorWal98}
{\sc Bordemann, M., Waldmann, S.:} {\it Formal GNS Construction and
States in Deformation Quantization.} Commun. Math. Phys. {\bf 195},
549--583 (1998).
%
\bibitem{Del95}
{\sc Deligne, P.:} {\it D{\'{e}}formations de l'Alg{\`{e}}bre des
       Fonctions d'une Vari{\'{e}}t{\'{e}} Symplectique:
       Comparaison entre Fedosov et DeWilde, Lecomte.} Sel. Math.,
       New Series {\bf 1 (4)}, 667--697 (1995).
%
\bibitem{DeWLec83b}
{\sc DeWilde, M., Lecomte, P. B. A.:} {\it Existence of
Star-Products and of Formal Deformations of the Poisson Lie Algebra
of Arbitrary Symplectic Manifolds.} Lett. Math. Phys. {\bf 7},
487--496 (1983).
%
\bibitem{Fed94}
{\sc Fedosov, B. V.:} {\it A Simple Geometrical Construction of
Deformation Quantization.} J. Diff. Geom. {\bf 40}, 213--238
(1994).
%
\bibitem{Fed96}
{\sc Fedosov, B. V.:} {\it Deformation Quantization and Index
Theory.} Akademie Verlag, Berlin (1996).
%
\bibitem{Fed98}
{\sc Fedosov, B. V.:} {\it Non-Abelian Reduction in Deformation
Quantization.} Lett. Math. Phys. {\bf 43}, 137--154 (1998).
%
\bibitem{GR99}
{\sc Gutt, S., Rawnsley, J.:} {\it Equivalence of star products on
a symplectic manifold; an introduction to Deligne's \v{C}ech
cohomology classes.} J. Geom. Phys. {\bf 29}, 347--392 (1999).
%
\bibitem{Gut00}
{\sc Gutt, S.:} {\it Variations on deformation quantization.} in:
{\sc Dito, G., Sternheimer, D. (eds.):} Conf\'{e}rence Mosh\'{e}
Flato 1999, Vol. I. Kluwer Academic Publ., Dordrecht, 217--254
(2000).
%
\bibitem{Hal99a}
{\sc Halbout, G.:} {\it Sur la classification des d\'{e}formations
des vari\'{e}t\'{e} de Poisson.} Th\`{e}se, Institut de recherche
math\'{e}matique avanc\'{e}e, Universit\'{e} Louis Pasteur et
C.N.R.S (1999).

\bibitem{Kar96}
{\sc Karabegov, A. V.:} {\it Deformation Quantization with
Separation of Variables on a K\"ahler Manifold.} Commun. Math.
Phys. {\bf 180}, 745--755 (1996).
%
\bibitem{Kar98}
{\sc Karabegov, A. V.:} {\it Cohomological Classification of
Deformation Quantization with Separation of Variables.} Lett. Math.
Phys. {\bf 43}, 347--357 (1998).
%
\bibitem{Kar99}
{\sc Karabegov, A. V.:} {\it On Fedosov's approach to Deformation
Quantization with Separation of Variables.} in: {\sc Dito, G.,
Sternheimer, D. (eds.):} Conf\'{e}rence Mosh\'{e} Flato 1999, Vol.
II. Kluwer Academic Publ., Dordrecht, 167--176 (2000).
%
\bibitem{KarSch00}
{\sc Karabegov, A. V., Schlichenmaier, M.:} {\it Identification of
Berezin-Toeplitz Deformation Quantization.} J. reine angew. Math.
{\bf 540}, 49--76 (2001).
%
\bibitem{Kas95}
{\sc Kassel, C.:} {\it Quantum Groups.} Graduate Texts in
Mathematics {\bf 155}, Springer Verlag, New York, Berlin,
Heidelberg (1995).
%
\bibitem{Kon97}
{\sc Kontsevich, M.:} {\it Deformation Quantization of Poisson
Manifolds, I.} Preprint, September 1997, {\bf q-alg/9709040}.
%
\bibitem{Neu99}
{\sc Neumaier, N.:} {\it Local $\nu$-Euler Derivations and
Deligne's Characteristic Class of Fedosov Star Products and Star
Products of Special Type.} Commun. Math. Phys. {\bf 230}, 271--288
(2002).
%
\bibitem{Neu02b}
{\sc Neumaier, N.:} {\it Some General Aspects of Fedosov-like
Products.} Preprint in preparation.
%
\bibitem{NeuWal02}
{\sc Neumaier, N.,Waldmann S.:} {\it Morita Equivalence Bimodules
for Wick Type Star Products.}  Preprint, July 2002, Freiburg
FR-THEP-2002/09, {\bf math.QA/0207162}. To appear in J. Geom. Phys.
%
\bibitem{NT95a} {\sc Nest, R., Tsygan, B.:} {\it Algebraic Index
Theorem.} Commun. Math. Phys. {\bf 172}, 223--262 (1995).
%
\bibitem{OmoMaeYos91}
{\sc Omori, H., Maeda, Y., Yoshioka, A.:} {\it Weyl Manifolds and
Deformation Quantization.} Adv. Math. {\bf 85}, 224--255 (1991).
%
\bibitem{Wal99}
{\sc Waldmann, S.:} {\it Zur Deformationsquantisierung in der
klassischen Mechanik: Observablen, Zust\"ande und Darstellungen.}
Ph.D. thesis, Fakult\"at f\"ur Physik der
Albert-Ludwigs-Universit\"at, Freiburg (1999). (available at:
http://idefix.physik.uni-freiburg.de/\verb|~|stefan/)
%
\bibitem{WX97}
{\sc Weinstein, A., Xu, P.:} {\it Hochschild cohomology and
characteristic classes for star-products.} in: {\sc Khovanskij, A.
et al. (eds.):} Geometry of differential equations. Dedicated to V.
I. Arnol'd on the occasion of his {$60^{\rm th}$} birthday.
Providence, Amer. Math. Soc. Transl., Ser. 2, {\bf 186 (39)},
177--194 (1998).
\end{thebibliography}
\end{document}